\documentclass[11pt]{article}
\usepackage{comment}
\usepackage[margin=1in]{geometry}
\usepackage[T1]{fontenc}
\usepackage{lmodern}
\usepackage{microtype}
\usepackage{amsmath,amssymb,amsfonts,amsthm,mathtools}
\usepackage{bm}
\usepackage{graphicx}
\usepackage{booktabs}
\usepackage{array}
\usepackage{enumitem}
\usepackage{natbib}
\usepackage{xcolor}
\usepackage{gensymb}
\usepackage{mathrsfs}
\usepackage[colorlinks=true,linkcolor=blue!55!black,citecolor=blue!55!black,urlcolor=blue!55!black]{hyperref}
\usepackage{cleveref}
\usepackage{authblk}
\newtheorem{theorem}{Theorem}[section]
\newtheorem{proposition}[theorem]{Proposition}
\newtheorem{lemma}[theorem]{Lemma}
\newtheorem{corollary}[theorem]{Corollary}

\theoremstyle{definition}

\theoremstyle{remark}
\newtheorem{remark}[theorem]{Remark}

\newcommand{\N}{\mathbb{N}}
\newcommand{\E}{\mathbb{E}}
\newcommand{\Pp}{\mathbb{P}}
\newcommand{\one}{\mathbf{1}}
\newcommand{\TV}{D_{\mathrm{TV}}}
\newcommand{\KL}{D_{\mathrm{KL}}}

\newcommand{\clip}{\operatorname{clip}}
\newcommand{\argmax}{\operatorname*{arg\,max}}
\newcommand{\argmin}{\operatorname*{arg\,min}}
\newcommand{\cP}{\mathcal{P}}
\newcommand{\cQ}{\mathcal{Q}}
\newcommand{\cE}{\mathcal{E}}
\newcommand{\cF}{\mathcal{F}}
\newcommand{\cX}{\mathcal{X}}
\newcommand{\cK}{\mathcal{K}}
\newcommand{\pos}[1]{\left(#1\right)_{+}}

\newcommand{\kl}{\operatorname{kl}}

\title{Beyond least-favorable distribution pairs:\\ asymptotically instance-optimal robust hypothesis testing}
\author[1]{Aytijhya Saha}
\author[2]{Aaditya Ramdas}

\affil[1]{Massachusetts Institute of Technology.  \texttt{aytijhya@mit.edu}}
\affil[2]{Stanford University. \texttt{aramdas@stanford.edu}}
\date{August 2026}

\begin{document}
\maketitle

\begin{abstract}
When testing a composite null $\cP$ against a composite alternative $\cQ$, it has been known since \cite{huber1973minimax} that if a least favorable distribution (LFD) pair exists, then a version of its likelihood ratio yields a minimax optimal test, and \cite{huber1965robust} derived a similar result in the sequential setting for a specific robust testing problem that is the focus of the current paper; the same likelihood ratio is log-optimal in the worst case; its expected logarithm under the alternative equals $\inf_{P \in \cP,S \in \cQ} \KL(S\|P)$, whose inverse also determines the stopping time of the test under the alternative. 
These results focus on optimizing worst-case performance, but the resulting methods have lower power if the data come from simpler alternatives. This paper studies an arguably more practical goal of achieving asymptotic \emph{instance optimality} for ``power-one'' sequential tests, which remarkably may even be possible when $\cQ = \cP^c$. This means that when the data are drawn from any $S \in \cP^c$, the test has power one, and its stopping time is determined by the larger information number $\inf_{P \in \cP} \KL(S\|P)$, which governs the performance of an oracle that knows $S$. To demonstrate this point, our paper focuses on the case of testing a TV ball around a discrete distribution with $d$ components, and designs a computationally and memory efficient algorithm that is also asymptotically instance-optimal. Specifically, we consider two different asymptotic settings: one is the “Wald setting'', where the
type-1 error level $\alpha$ approaches zero for a fixed alternative $S$ and the second is
the “Farrell setting”, where $\alpha$ is fixed and $\inf_{P \in \cP} \KL(S\|P)$ approaches $0$. Our sequential test is simultaneously optimal in both regimes.
\end{abstract}

\section{Introduction}

Consider sequential testing of a composite null class $\cP$ of probability laws on a measurable space $\cX$.  In the classical formulation, one specifies a composite alternative $\cQ$ and seeks a test with strong worst-case performance over $\cQ$.  If there is a least-favorable pair $(P^\star,Q^\star)\in\cP\times\cQ$, then the corresponding likelihood ratio provides a minimax test: it is designed against the hardest alternative $Q^\star$ and protects against the hardest null $P^\star$.  Robust likelihood-ratio tests of this form go back to Huber's censored probability-ratio construction for neighborhood models~\citep{huber1965robust}, while Huber and Strassen later developed a broad minimax theory in terms of least-favorable pairs and capacities~\citep{huber1973minimax}.  More recently, \cite{saha2025huber} formulated this principle in terms of e-values for arbitrary composite nulls and alternatives and showed how one can perform one-sided sequential tests for Huber-contamination and total-variation neighborhoods, including sequentially adaptive contamination.

The minimax criterion is compelling when performance at the hardest alternative is the primary objective.  It can, however, leave substantial power on the table when the realized data-generating distribution is easier than the least-favorable one.  A test optimized for $Q^\star$ has no reason to exploit the additional separation of a particular $Q\in\cQ$ from the null class.  In sequential problems, this loss appears directly as slower exponential growth of the evidence process and hence a larger detection delay.  Moreover, in many testing problems there is no naturally separated alternative class for which a nontrivial least-favorable pair exists.  A canonical example is testing a composite null against its complement.

This motivates a different objective.  Instead of fixing a composite class and optimizing the worst case over it, take the alternative to be simply the complement
\[
    \cQ=\cP^c,
\]
and ask for a single valid sequential test that adapts to whichever $S\in\cP^c$ actually generates the data, given a level $\alpha \in (0,1)$. We construct an e-process and a ``power-one" sequential test that can control the type-1 error at a prespecified level $\alpha$ and, for a fixed alternative $S\notin\cP$, the e-process asymptotically achieves the oracle growth rate, 
\begin{equation}\label{eq:intro-information-general}
    \Gamma(S)
    :=
    \inf_{R\in\cP}\KL(S\|R),
\end{equation}
which is known to be optimal for testing $\cP$ \citep{lardy2024reverse,larsson2025numeraire}. Here $\KL$ denotes the KL divergence. See \Cref{sec:background} for a detailed discussion on the e-process and its growth rate.

We develop this idea for a canonical robust-testing model. For some $d\in\N$, let $\cX=\{1,\ldots,d\}$ be the observation space and consider a discrete distribution with $d$ components with probability mass function
\[
    P=(P_1,\ldots,P_d),
\]
with $P_i>0,$ for all $i=1,\cdots,d$.
Take the null hypothesis to be the total-variation neighborhood \citep{huber1965robust} around $P$:
\begin{equation}\label{eq:intro-tv-null}
    \cP_0
    :=
    B_{\mathrm{TV}}(P,\epsilon_0)
    =
    \{R\in\Delta_d:\TV(R,P)\leq\epsilon_0\},
\end{equation}
where $\Delta_d$ is the probability simplex, $\TV(R,P)$ is the total variation distance between $R$ and $P$:
\[
    \TV(R,P)=\frac12\sum_{i=1}^d|R_i-P_i|.
\]
With a slight abuse of notation, we will denote the probability distributions by their $d$-dimensional probability mass functions. The radius $\epsilon_0$ has a direct interpretation: under any distribution in the null, the probability of every event differs from its nominal probability by at most $\epsilon_0$. It therefore specifies the discrepancies against which the test must protect, even when evidence is monitored continuously.

But in contrast to \cite{huber1965robust,saha2025huber}, our alternative is not another prespecified ball; it is the entire complement $\Delta_d\setminus\cP_0$.
Our goal here is to approach the instance-optimal information rate defined in \eqref{eq:intro-information-general}, instead of the minimax information rate, $\inf_{P \in \cP,S \in \cQ} \KL(S\|P)$ as in \cite{saha2025huber}.

We also construct e-processes for both Wald and Farrell optimality for testing
$\cP_0$ against its entire complement
$\cP_0^c=\Delta_d\setminus\cP_0$. Specifically, we seek a single e-process whose threshold-crossing times define
level-$\alpha$ tests with power one under every alternative, and whose
expected stopping times are optimal in two complementary regimes.
In the Wald regime, the alternative $S\notin\cP_0$ is fixed and
$\alpha\downarrow0$; for every fixed full-support $S$, we aim to attain
\[
    \mathbb{E}_S[\tau_\alpha]
    \sim \frac{\log(1/\alpha)}{\Gamma(S)}.
\]
In the Farrell regime, $\alpha\in(0,1/2)$ is fixed and the alternative
approaches the null boundary, so that $\Gamma(S)\downarrow0$; with
$P,\epsilon_0,d$ fixed, we aim for
\[
    \mathbb{E}_S[\tau_\alpha]
    =
    O\!\left(
        \frac{\log\log(1/\Gamma(S))}{\Gamma(S)}
    \right).
\]
The first guarantee attains the exact information-theoretic leading
constant, while the second matches the Farrell lower bound in order
in its sequence-wise sense. We obtain both guarantees by averaging
an adaptive clipped-Kelly e-process with a complementary e-process
designed for alternatives near the null boundary.

The remainder of this introduction reviews the required background, discusses related work, and states these contributions precisely.

\subsection{Background on e-values, test supermartingales, one-sided tests, growth rate and consistency}
\label{sec:background}
An e-variable is a nonnegative random variable whose expectation is bounded above by one under every distribution in the null class. Specifically, $B$ is an e-variable for $\mathcal P$ if
\[
    \mathbb E_P[B]\leq1
    \qquad\text{for all }P\in \mathcal P.
\]
The observed realization of an e-variable is called an e-value.

Let $M=\{M_n\}_{n\geq0}$ be an integrable process adapted to a filtration $\mathscr F=\{\mathscr F_n\}_{n\geq0}$. We call $M$ a \emph{martingale} for $\mathcal P$ if, for every $n\geq1$,
\begin{equation}\label{eq:martingale}
    \mathbb E_P[M_n\mid\mathscr F_{n-1}]
    =
    M_{n-1},
    \qquad \forall P\in \mathcal P.
\end{equation}
Similarly, $M$ is a \emph{supermartingale} for $\mathcal P$ if, for every $n\geq1$,
\begin{equation}\label{eq:supermartingale}
    \mathbb E_P[M_n\mid\mathscr F_{n-1}]
    \leq
    M_{n-1},
    \qquad \forall P\in \mathcal P.
\end{equation}
A \emph{test (super)martingale} is a process that satisfies the corresponding condition under every distribution in $\mathcal P$, is nonnegative, and starts at $M_0=1$. A stopping time $\tau$ is a nonnegative integer-valued, possibly infinite random variable satisfying $\{\tau\leq n\}\in\mathscr F_n$ for every $n\in\mathbb N$. We write $\mathcal T$ for the collection of all such stopping times, including those that may never stop.

Ville's inequality~\citep{ville1939etude} guarantees that a test (super)martingale satisfies
\begin{equation}
    \sup_{P\in \mathcal P}
    P\!\left(\exists n\in\mathbb N:M_n\geq1/\alpha\right)
    \leq\alpha.
\end{equation}
A short proof appears in \citet[Lemma~1]{howard2020time}. Equivalently,
$P(M_\tau\geq1/\alpha)\leq\alpha$ for every $\tau\in\mathcal T$ and every $P\in \mathcal P$; see \citet[Lemma~3]{howard2021time}. Consequently, stopping and rejecting at
\begin{equation}\label{eq:stopping}
    \tau_\alpha
    =
    \inf\{n\geq1:M_n\geq1/\alpha\}
\end{equation}
defines a level-$\alpha$ sequential test. Under the null, the probability of ever stopping with a false rejection is at most $\alpha$; equivalently, $\tau_\alpha=\infty$ with probability at least $1-\alpha$.

This formulation differs from Wald's classical sequential testing framework~\citep{wald1945sequential}, in which the procedure may terminate by either accepting or rejecting the null, with the stopping rule calibrated to specified Type-I and Type-II error probabilities. Our setting instead follows Robbins' theory of ``power-one'' or one-sided tests~\citep{darling1968some,robbins1970statistical}: only a target Type-I error level is specified, and sampling stops exclusively to reject the null. Such a test is called consistent under an alternative $Q$ if it has power one under $Q$.

A test (super)martingale is called consistent, or power one, for an alternative $Q$ if
\[
    M_n\longrightarrow\infty
    \qquad Q\text{-almost surely}.
\]
Under each alternative, the process then accumulates unbounded evidence against the null and eventually crosses every finite rejection threshold. In i.i.d.\ settings, evidence can accumulate exponentially, as emphasized by \citet{kelly1956new}. This motivates measuring performance through the worst-case logarithmic growth rate
\[
    \liminf_{n\to\infty}\frac{\log M_n}{n}
\]
\citep{shafer2021testing,grunwald2024safe,waudby2023estimating}. A positive growth rate implies consistency. A test (super)martingale attaining the largest possible growth rate is called log-optimal. This optimality criterion has been studied and advocated by \citet{breiman1961optimal}, \citet{cover1987log}, \citet{shafer2021testing}, and \citet{grunwald2024safe}, among others.

\cite{agrawal2025stopping} provides upper and lower bounds for power-one sequential tests in two different asymptotic settings: one is the “Wald setting'', where the
type-1 error level $\alpha$ approaches zero for a fixed alternative $S$. The second lower bound applies to
the “Farrell setting”, where $\alpha$ is fixed and $\Gamma(S)$ approaches 0 along a sequence of alternatives.


\subsection{Related work}

\paragraph{Robust/minimax hypothesis testing.}
Robust likelihood-ratio testing begins with Huber's censored probability-ratio test~\citep{huber1965robust}, which clips individual likelihood-ratio factors to prevent a small number of contaminated observations from dominating the evidence.  \cite{huber1973minimax} developed a general minimax theory of hypothesis testing under capacities and least-favorable distributions.  These works establish the classical principle that, when an appropriate least-favorable pair exists, testing can be reduced to a likelihood-ratio test between that pair.

\paragraph{E-values and e-processes.}
The second thread of literature that is relevant involves recent fundamental advances in composite hypothesis testing. For example, 
 \cite{grunwald2024safe} introduce growth-rate optimal (GRO) e-variables
for testing problems with composite null and alternative. The universal inference method \citep{wasserman2020universal} used an \emph{e-value} to propose a level-$\alpha$ test for any composite null without requiring any regularity conditions, but requiring reference measures. 
 Reverse information projections and numeraire e-variables characterize pointwise log-optimal evidence under composite nulls~\citep{lardy2024reverse,larsson2025numeraire}. 
 See \cite{ramdas2022game} for a detailed review.

 Our contribution is to turn this static pointwise characterization into an adaptive sequential procedure with explicit stopping-time and computational guarantees in the finite-alphabet robust model.
 The closest modern precursor is \cite{saha2025huber}.  They prove that least-favorable pairs yield log-optimal e-values for arbitrary composite nulls and alternatives, construct Huber-robust and TV-robust versions of such tests, and obtain sequential validity even when the contamination mechanism varies adaptively with the past, extending~\cite{huber1973minimax}.
 In particular, their theory already gives the least-favorable clipped likelihood ratio for two separated TV neighborhoods.  Our question is orthogonal: rather than optimizing the worst case over a prespecified alternative class, we take the alternative to be $\cP^c$ and ask whether one valid procedure can asymptotically attain the oracle rate \eqref{eq:intro-information-general} at every realized alternative.

 \cite{ram2026power} show the existence of relatively growth-rate optimal e-processes and thus power-one tests under rather general conditions, including the TV ball. However, their work is abstract and their constructions are far from practical, and they do not have stopping time guarantees or regret bounds. Our work extends theirs by adding the specific contributions (i), (ii), (iii), and (iv) mentioned in the next subsection.

 \subsection{Contributions}

\begin{enumerate}[label=(\roman*)]

\item \textbf{Finite-sample adaptive regret.}
For testing the null $\mathcal P$ in \eqref{eq:intro-tv-null} against an unknown full-support alternative $S\in\mathcal P^c$, we construct an e-process $\{M_n^{\mathrm{CK}}\}_n$ and show the finite-sample bound on the regret
\begin{equation}\label{eq:intro-plugin-regret}
    0\leq
    n\Gamma(S)-\E_S\log M_n^{\mathrm{CK}}
    \leq
    \log\binom{n+d-1}{d-1}.
\end{equation}
Consequently, the expected regret is $O(d\log n)$ and the normalized expected and almost-sure log wealth converge to $\Gamma(S)$.  

\item \textbf{Wald-optimal expected stopping time.}
For the level-$\alpha$ test using the that e-process,
$\tau_\alpha^{\mathrm{CK}}
    :=
    \inf\{n:M_n^{\mathrm{CK}}\geq\alpha^{-1}\},$
we prove that
\[
    \E_S\tau_\alpha^{\mathrm{CK}}
    =
    \frac{\log(1/\alpha)}{\Gamma(S)}(1+o(1)),
    \qquad \alpha\downarrow0.
\]
It is well known that no level-$\alpha$ sequential test can improve the first-order constant $1/\Gamma(S)$ at the same instance \citep{agrawal2025stopping}.

\item \textbf{Farrell-optimal expected stopping time.}
We construct another e-process $\{M^F_n\}_n$ for adaptation near the null boundary
for which the level $\alpha$ stopping time, $\tau_\alpha^{\mathrm{F}}$, is Farrell-optimal:
\[
\E_S\tau_\alpha^{\mathrm{F}}= O_{P,\epsilon_0,d}\!\left(
 \frac{d+\log(1/\alpha)+\log\log(e^e/\Gamma(S))}{\Gamma(S)}
 \right),\qquad 0<\Gamma(S)\le1,
\]
that is, for fixed $\alpha\in(0,1/2),$  we have $\E_S\tau_\alpha^{\mathrm{F}}=O(\Gamma(S)^{-1}\log\log(1/\Gamma(S)))$, as  $\Gamma(S)\downarrow 0$.
 At fixed $\alpha\in(0,1/2)$, it matches the optimal Farrell order in the sequence-wise sense of the lower bound. Averaging this process with clipped Kelly retains the exact Wald constant and increases the expected-regret bound by at most $\log2$.

\item \textbf{Exact logarithmic-time implementation.}
A naive fresh oracle computation sorts all $d$ ratios and costs $O(d\log d)$ per round.  We show that, after a rescaling, only the coordinate corresponding to the newly observed symbol changes.  Maintaining these quantities in an augmented balanced binary search tree yields the \emph{same} fresh every-round clipped-Kelly oracle in $O(\log d)$ worst-case update time, after $O(d\log d)$ initialization, using $O(d)$ memory.
\end{enumerate}

  Our adaptive process instead learns the oracle corresponding to the realized alternative, enjoys logarithmic finite-sample regret, achieves the information-theoretically optimal first-order stopping time at that instance, and can be implemented efficiently.

\paragraph{Organization.}
\Cref{sec:adaptive} constructs an adaptive e-process showing that its expected stopping time achieves Wald optimality as the type-I error converges to zero. \Cref{sec:farrell} studies adaptation near the null boundary and Farrell-type optimality, and combines the two regimes.  \Cref{sec:dynamic} gives efficient dynamic implementation of both e-processes. \Cref{sec:expt} shows a simulation experiment, and \Cref{sec:discussion} discusses the scope of the results. Proofs of the theoretical results  appear in \Cref{app:proofs}. An alternative
oracle-mixture construction is provided in \Cref{app:mixture}.

\section{Adaptive clipped-Kelly e-process}
\label{sec:adaptive}

From this point onward we specialize to the finite-alphabet TV-robust model in \eqref{eq:intro-tv-null}.  We assume throughout that $d\geq2$, $P_i>0$ for every $i$, and
\[
    0<\epsilon_0<1-\min_i P_i,
\]
so that the TV ball is nontrivial and its complement is nonempty.  The case $\epsilon_0=0$ reduces to the ordinary likelihood-ratio/Kelly construction and is simpler than the robust case considered here.
\subsection{The pointwise Kelly oracle}
\label{sec:oracle}

  We write
\[
    \Delta_d^\circ:=\{U\in\Delta_d:U_i>0\text{ for all }i\}
\]
for the interior of the simplex.  All oracle parameters used by the adaptive procedures lie in $\Delta_d^\circ$ because the predictors use strictly positive pseudocounts.

Define the one-step e-variable class
\begin{equation}\label{eq:E0}
    \cE_0
    :=\left\{E\geq 0:
      \sup_{R\in\cP_0}\E_R E\leq1\right\}.
\end{equation}

Note that for $U\in\Delta_d$, $\Gamma(U)$ defined in \eqref{eq:intro-information-general} can be written as
\begin{equation}\label{eq:Gamma-def}
    \Gamma(U)
    =\sup_{E\in\cE_0}\sum_{i=1}^d U_i\log E_i.
\end{equation}

To construct a sequential test that adapts to an unknown alternative, we first
solve the idealized problem in which the alternative distribution $U$ is known.
An oracle with this knowledge would choose a betting factor that is valid for
every distribution in the null class $\cP_0$ and maximizes the expected
logarithmic return under $U$:
\[
    E_U \in \arg\max_{E\in\mathcal E_0}
    \sum_{i=1}^d U_i\log E_i.
\]
Under i.i.d.\ sampling from $U$, repeatedly using this factor yields expected
log wealth $n\Gamma(U)$ after $n$ observations. Thus the pointwise oracle
identifies both the information rate available at a particular alternative
and the betting factor that an adaptive procedure should aim to learn.

The main result of this section gives an explicit characterization of this
oracle for a TV-ball null. We show that the optimal factor is obtained by
clipping the nominal likelihood ratio $U_i/P_i$ between two thresholds
determined by the robustness radius. Equivalently, it is the likelihood ratio
$U_i/R_U(i)$, where $R_U$ is the reverse-KL projection of $U$ onto $\cP_0$.
This characterization provides the link between the static optimization
problem and the adaptive construction below: we can apply the oracle map
to an estimate of $U$ based on past observations, while retaining null
validity for every such estimate.

We work with strictly positive oracle inputs $U\in\Delta_d^\circ$, as is
automatically the case for the smoothed empirical predictors used later.
If $U\in\cP_0$, no null-valid factor can have positive expected log return
under $U$, and the constant factor one is optimal. Accordingly, for
$U\in\Delta_d^\circ\cap\cP_0$, set $E_U\equiv1$, $R_U=U$, and $\Gamma(U)=0$.
For the remainder of this section, fix
$U\in\Delta_d^\circ\setminus\cP_0$. Strict positivity also ensures that the
lower clipping threshold introduced below is strictly positive.

Define $a(U)\in(0,1)$ and $b(U)\in(1,\infty)$ through
\begin{align}
    \sum_{i=1}^d\pos{a(U)P_i-U_i}
    &=a(U)\epsilon_0,
    \label{eq:oracle-low}\\
    \sum_{i=1}^d\pos{U_i-b(U)P_i}
    &=b(U)\epsilon_0.
    \label{eq:oracle-high}
\end{align}
Set
\begin{equation}\label{eq:oracle-e}
    E_U(i)=\clip\!\left(\frac{U_i}{P_i},a(U),b(U)\right),
    \qquad
    \clip(z,a,b)=\min\{b,\max\{a,z\}\},
\end{equation}
and
\begin{equation}\label{eq:reverse-projection}
    R_U(i)=\frac{U_i}{E_U(i)}.
\end{equation}

\begin{lemma}[Existence and uniqueness of the oracle thresholds]\label{lem:oracle-roots}
If $U\in\Delta_d^\circ\setminus\cP_0$, then \eqref{eq:oracle-low}--\eqref{eq:oracle-high} have unique solutions $a(U)\in(0,1)$ and $b(U)\in(1,\infty)$.
\end{lemma}

\begin{theorem}[Reverse projection and pointwise Kelly oracle]\label{thm:oracle}
Let $U\in\Delta_d^\circ$.  Then:
\begin{enumerate}[label=(\alph*)]
\item $R_U\in\cP_0$ and $E_U\in\cE_0$;
\item $R_U$ is the unique reverse information projection (RIPr) and
$ R_U=\argmin_{R\in\cP_0}\KL(U\|R);$
\item
$\Gamma(U)=\KL(U\|R_U)
    =\inf_{R\in\cP_0}\KL(U\|R),$
and $E_U=U/R_U$ is the unique maximizer in \eqref{eq:Gamma-def}.
\end{enumerate}
\end{theorem}

These facts specialize the reverse-projection characterization of log-optimal e-variables \cite{lardy2024reverse,larsson2025numeraire}. In particular, it turns out that in this case we can avoid the explicit use of their more general bipolar machinery, because we are able to find an e-variable which is the likelihood ratio of $U$ to some element of $\cP_0$, immediately making it the log-optimal e-variable by results in the latter paper.

\begin{corollary}\label{cor:ftrl}
For positive weights $c_1,\ldots,c_d$ and
$U_i=c_i/\sum_jc_j$,
\[
    E_U=\argmax_{E\in\cE_0}\sum_{i=1}^d c_i\log E_i.
\]
\end{corollary}

The pointwise oracle identifies the best expected logarithmic growth
attainable when the alternative distribution is known. But the oracle is unavailable because the alternative is unknown. We now construct
a single sequential procedure that learns this oracle from the observations.
At each round, we estimate the alternative using past data and apply the
oracle map to that estimate. Since $E_U$ is valid for the entire null class
for every oracle input $U$, this predictable substitution preserves
anytime validity. The remaining task is to quantify the loss in log wealth
caused by having to estimate the alternative.

For the expected-regret and stopping-time analysis, we assume that the observations are i.i.d.\
from a fixed alternative $S$. For $S\in\Delta_d$ and
$U\in\Delta_d^\circ$, write
\[
    h(S,U):=\sum_{i=1}^d S_i\log E_U(i).
\]
This is the expected logarithmic return under $S$ of the betting factor
designed for $U$. In particular, $h(S,S)=\Gamma(S)$ when $S$ has full
support, so the difference
\[
    \Gamma(S)-h(S,U)
\]
measures the one-step cost of using $U$ in place of the true alternative.

The following lemma bounds this cost by $D_{\mathrm{KL}}(S\|U)$.
It provides the key link between learning the alternative and accumulating
evidence: a predictor with small cumulative log-loss regret also yields
a robust e-process with small expected regret relative to the oracle.
Below we obtain a stronger pathwise finite-sample bound for a smoothed
empirical predictor.

\begin{lemma}[Oracle mismatch is bounded by KL]\label{lem:oracle-mismatch}
Let $S\notin\cP_0$ have full support and let $U\in\Delta_d$ have strictly positive coordinates.  Then
\begin{equation}\label{eq:oracle-mismatch}
    0
    \leq
    \Gamma(S)-h(S,U)
    \leq
    \KL(S\|U).
\end{equation}
\end{lemma}
\subsection{Predictable clipped-Kelly process}

Fix pseudocounts $\alpha_i>0$ and let
\[
    \alpha_+=\sum_i\alpha_i,
    \qquad
    N_i(t)=\sum_{r=1}^t\one\{X_r=i\}.
\]
The every-round predictor and factor are
\begin{equation}\label{eq:plugin}
    \widehat U_{t-1,i}
    =\frac{N_i(t-1)+\alpha_i}{t-1+\alpha_+},
    \qquad
    E_t=E_{\widehat U_{t-1}},
\end{equation}
and
\[
    M_n^{\mathrm{CK}}=\prod_{t=1}^n E_t(X_t).
\]
By Corollary~\ref{cor:ftrl}, this is exactly a follow-the-regularized-leader
or Kelly update over $\cE_0$.

The fresh predictor appears to require an $O(d\log d)$ oracle recomputation after every observation.  Section~\ref{sec:dynamic} shows that the identical every-round factor can instead be maintained exactly in $O(\log d)$ time.

\begin{theorem}
\label{thm:adaptive-valid}
The process $(M_n^{\mathrm{CK}})_n$ is a nonnegative test supermartingale under $\cP$.
\end{theorem}

This validity statement holds for every choice of positive
pseudocounts and requires no accuracy guarantee for the predictor.
Conditional on the past, the selected factor belongs to
$\mathcal E_0$, which is sufficient for the supermartingale property.
The quality of the predictor instead determines how quickly evidence
accumulates under an alternative. We next quantify this learning cost.

For a particularly clean finite-sample bound, we now take the Laplace
pseudocounts
\begin{equation}\label{eq:laplace-prior}
    \alpha_1=\cdots=\alpha_d=1.
\end{equation}
Other fixed positive pseudocounts lead to the same logarithmic order with
constants depending on the pseudocount vector.

\begin{theorem}[Pathwise regret of adaptive clipped Kelly]\label{thm:plugin-regret}
Assume \eqref{eq:laplace-prior}.  For every observation sequence
$x_1,x_2,\ldots\in\cX$ and every $n\geq1$, let
$\bar S_{n,i}=N_i(n)/n$ and define the best fixed null-valid wealth in hindsight by
\[
    W_n^\star
    :=\max_{E\in\cE_0}\prod_{t=1}^n E(x_t).
\]
Then
\[
    0\leq
    \log W_n^\star-\log M_n^{\mathrm{CK}}
    =n\Gamma(\bar S_n)-\log M_n^{\mathrm{CK}}
    \leq\log\binom{n+d-1}{d-1}.
\]

\end{theorem}
The benchmark $W_n^\star$ is allowed to inspect all $n$ observations
before choosing a single null-valid factor to use throughout the
sequence. The theorem shows that the predictable procedure nearly
matches this hindsight choice: for fixed $d$, the difference in
log wealth grows only logarithmically with the horizon. The comparison
holds for every observation sequence.

The regret penalty depends only on $n$ and $d$; the robustness radius
enters through the feasible factors and the benchmark
$\Gamma(\bar S_n)$. Under i.i.d.\ sampling from $S$, the fixed oracle
$E_S$ is one of the competitors defining $W_n^\star$. Taking
expectations therefore gives the upper bound in the following
corollary.

\begin{corollary}[Expected regret and pointwise instance optimality]\label{cor:plugin-expected-growth}
Assume \eqref{eq:laplace-prior} and that the observations are i.i.d.\ from a
fixed full-support $S\in\Delta_d\setminus\cP_0$.  For every $n\geq1$,
\begin{equation}\label{eq:plugin-regret}
    0
    \leq
    n\Gamma(S)-\E_S\log M_n^{\mathrm{CK}}
    \leq
    \log\binom{n+d-1}{d-1}.
\end{equation}
Consequently, for $d\geq2$,
\begin{equation}\label{eq:plugin-regret-simple}
    n\Gamma(S)-\E_S\log M_n^{\mathrm{CK}}
    \leq
    (d-1)\log\!\left(\frac{e(n+d-1)}{d-1}\right).
\end{equation}
In particular,
\begin{equation}\label{eq:plugin-expected-growth}
    0
    \leq
    \Gamma(S)-\frac1n\E_S\log M_n^{\mathrm{CK}}
    \leq
    \frac1n\log\binom{n+d-1}{d-1}
    \longrightarrow0.
\end{equation}
\end{corollary}

\begin{corollary}[Almost-sure pointwise instance optimality]\label{cor:plugin-as-growth}
Under the assumptions of Corollary~\ref{cor:plugin-expected-growth},
\begin{equation}\label{eq:plugin-as-growth}
    \frac1n\log M_n^{\mathrm{CK}}
    \longrightarrow\Gamma(S)
    \qquad\text{almost surely under }S.
\end{equation}
\end{corollary}
Together, these corollaries show that the cost of learning the
alternative vanishes per observation, both in expectation and almost
surely. Since $\Gamma(S)>0$ for every fixed alternative under
consideration, the wealth diverges almost surely and the test
eventually rejects at every fixed significance level. Thus the same
procedure has power one at every full-support alternative outside
the null ball. We next ask how many observations this eventual
rejection requires.

\subsection{Instance-optimal expected stopping time}\label{sec:stopping}

For an anytime-valid wealth process $M$, define
\begin{equation}\label{eq:tau}
    \tau_\alpha(M)
    :=\inf\{n\geq1:M_n\geq\alpha^{-1}\}.
\end{equation}
This gives a level-$\alpha$ test by Ville's inequality.
The growth benchmark $\Gamma(S)$ has a direct operational interpretation:
its reciprocal is the optimal first-order expected sample size per unit of
log evidence.

\begin{theorem}[Universal information lower bound]\label{thm:stop-lower}
Fix a full-support $S\notin\cP_0$, let $R_S$ be its reverse projection from
Theorem~\ref{thm:oracle}, and let $\tau$ be any stopping rule satisfying
\[
\sup_{R\in\cP_0}\Pp_R(\tau<\infty)\leq\alpha,
    \qquad
    \Pp_S(\tau<\infty)=1.
\]
Then
\begin{equation}\label{eq:stop-lower}
    \E_S\tau
    \geq\frac{\log(1/\alpha)}{\Gamma(S)}.
\end{equation}
\end{theorem}

The above information lower bound is a specialization of Theorem 3.1 of \cite{agrawal2025stopping}. We provide a self-contained proof for our special case in the appendix.

\begin{theorem}[Expected stopping time of adaptive clipped Kelly]\label{thm:stop-upper}
Assume the pseudocounts \eqref{eq:laplace-prior}.  Fix a full-support $S\notin\cP_0$ and let
\[
    \tau_\alpha^{\mathrm{CK}}
    :=\tau_\alpha(M^{\mathrm{CK}}).
\]
Then, as $\alpha\downarrow0$,
\begin{equation}\label{eq:stop-asymptotic}
    \E_S\tau_\alpha^{\mathrm{CK}}
    \sim
    \frac{\log(1/\alpha)}{\Gamma(S)}.
\end{equation}
\end{theorem}

\begin{remark}
The convergence $n^{-1}\E_S\log M^{\mathrm{CK}}_n\to\Gamma(S)$ does not by itself imply
\eqref{eq:stop-asymptotic}.  The proof of the theorem in the appendix uses the stronger exponential
concentration available on a finite alphabet to control the tail of the
hitting time and hence establish uniform integrability at the relevant scale.
\end{remark}

\section{Farrell optimal e-process: adaptation near the null boundary}
\label{sec:farrell}

The stopping-time result in the previous section fixes an alternative $S$ and lets the error level $\alpha$ decrease to zero. A different difficulty arises when $\alpha$ is fixed and the alternative approaches the null. In this regime, $\Gamma(S)$ decreases to zero, and the procedure must learn both where the null is violated and the scale of that violation. Optimality in the first limit does not imply optimality in the second. We now construct a complementary process that adapts to this unknown scale. Averaging it with clipped Kelly gives one test with the exact Wald constant and the optimal Farrell order in their respective limits.

Retain the standing assumptions on $P,d,\epsilon_0$, which are fixed throughout the small-gap analysis. Put
\[
 \Delta(S)=\TV(S,P)-\epsilon_0,\qquad
 \Gamma(S)=\inf_{R\in\cP_0}\KL(S\|R),\qquad
 \ell(x)=\log\log(e^e/x)\quad(0<x\le1).
\]
Performance is evaluated under i.i.d. observations from $S\notin\cP_0$. The small-gap bounds allow zero coordinates in $S$.

\subsection{The geometry of a small violation}
For an event $A\subseteq\cX$, every null distribution satisfies
\[
 R(A)\le P(A)+\epsilon_0.
\]
Moreover, $\TV(S,P)=\max_A\{S(A)-P(A)\}$. Thus an alternative has an event whose null constraint is violated by exactly $\Delta(S)$. The following lemma relates this event-level violation to the information number used in our earlier results.

\begin{lemma}[Uniform information comparison]
\label{lem:farrell-geometry}
Define
\[
 \eta_0=\min_{A\subseteq\cX:\,P(A)>\epsilon_0}
             \{P(A)-\epsilon_0\},\qquad
 c_0=\min\{\epsilon_0,\eta_0\}>0.
\]
For every $S\in\Delta_d\setminus\cP_0$,
\begin{equation}
 2\Delta(S)^2\le\Gamma(S)\le c_0^{-2}\Delta(S)^2.
 \label{eq:farrell-geometry}
\end{equation}
\end{lemma}

These constants depend on the fixed null geometry and need not be computed by the algorithm. In particular, the lemma does not assume a positive lower bound on the coordinates of $S$. It does not assert uniform constants as $d$ or $\epsilon_0$ vary. 

\subsection{Learn an event, bet for one block, and retain the capital}
Let $m_j=2^j$ for $j\ge0$, and let $\bar S_m$ be the empirical distribution after $m$ observations. At time $m_j$, select
\begin{equation}
 \begin{split}
 \widehat A_j&=\{i:\bar S_{m_j,i}>P_i\},\qquad q_j=P(\widehat A_j)+\epsilon_0,\\
 \widehat\Delta_j&=(\TV(\bar S_{m_j},P)-\epsilon_0)_+,
 \qquad \lambda_j=4\widehat\Delta_j.
 \end{split}
 \label{eq:farrell-witness}
\end{equation}
The selected event maximizes the empirical violation of a null constraint. Its betting intensity is proportional to the estimated excess. If the estimated excess is zero, the account does not bet.

Account $j$ starts at one, bets only during observations $m_j+1,\ldots,2m_j$, and then freezes. Precisely, let $\nu_j(n)=((n\wedge2m_j)-m_j)_+$ and set
\begin{equation}
 W_n^{(j)}=
 \begin{cases}
 1,&n\le m_j,\\
 \displaystyle\exp\left\{
 \lambda_j\sum_{t=m_j+1}^{n\wedge2m_j}
       (\one\{X_t\in \widehat A_j\}-q_j)
       -\frac{\nu_j(n)\lambda_j^2}{8}\right\},&n>m_j.
 \end{cases}
 \label{eq:farrell-account}
\end{equation}
Allocate initial capital $w_j=1/[(j+1)(j+2)]$ and define
\begin{equation}
 M_n^{\mathrm{F}}=\sum_{j=0}^\infty w_jW_n^{(j)}.
 \label{eq:farrell-process}
\end{equation}
The weights sum to one, so $M_0^{\mathrm{F}}=1$. The process is independent of $\alpha$.

At each time only one account is active. Learning the event from the preceding pilot keeps its choice predictable. Under an alternative, a sufficiently accurate pilot identifies an event with true violation of order $\Delta(S)$, and the next block then earns log wealth of order $m_j\Delta(S)^2$. Trying geometrically increasing pilot sizes costs only $\log(1/w_j)=O(\log(j+1))$ log capital. This becomes the iterated logarithm in the stopping-time bound.

\begin{theorem}[Small-gap sample complexity]
\label{thm:farrell-upper}
The process $M^{\mathrm{F}}$ is a test supermartingale under $\cP_0$. There is an absolute constant $C$ such that, for every $S\notin\cP_0$ and $\alpha\in(0,1)$,
\begin{equation}
 \E_S\tau_\alpha(M^{\mathrm{F}})
 \le\frac{C}{\Delta(S)^2}
   \left[d+\log(1/\alpha)+\ell(\Delta(S)^2)\right].
 \label{eq:farrell-delta-bound}
\end{equation}
Consequently, for $0<\Gamma(S)\le1$,
\begin{equation}
 \E_S\tau_\alpha(M^{\mathrm{F}})
 \le\frac{C_{P,\epsilon_0,d}}{\Gamma(S)}
   \left[d+\log(1/\alpha)+\ell(\Gamma(S))\right],
 \label{eq:farrell-gamma-bound}
\end{equation}
where $C_{P,\epsilon_0,d}<\infty$ depends only on the fixed null geometry. Also, $M_n^{\mathrm{F}}\to\infty$ almost surely under every such $S$.
\end{theorem}

For fixed $\alpha$, the theorem gives $O(\Gamma(S)^{-1}\log\log(1/\Gamma(S)))$ expected sample size uniformly as $\Gamma(S)\downarrow0$. The proof controls the entire hitting-time tail, rather than only the expected log wealth. No small-gap optimality claim for the original clipped-Kelly process is needed.

\subsection{One process for both asymptotic regimes}
Define
\begin{equation}
 M_n^{\mathrm{mix}}=\tfrac12M_n^{\mathrm{CK}}+\tfrac12M_n^{\mathrm{F}}.
 \label{eq:farrell-combined}
\end{equation}
Here $M^{\mathrm{CK}}$ is the Laplace-smoothed clipped-Kelly process defined in \eqref{eq:plugin}. The mixture is a simple device for retaining each component's guarantee; the same general combination principle appears in \cite[Proposition 5.1]{agrawal2025stopping}.

\begin{corollary}[Exact Wald constant and optimal Farrell order]
\label{cor:farrell-joint}
The process $M^{\mathrm{mix}}$ is an $\alpha$-independent test supermartingale under $\mathcal P_0$ and is power one under every i.i.d. $S\notin\cP_0$. For each fixed full-support $S\notin\cP_0$,
\begin{equation}
 \E_S\tau_\alpha(M^{\mathrm{mix}})\sim\frac{\log(1/\alpha)}{\Gamma(S)}
 \quad\text{as }\alpha\downarrow0.
 \label{eq:farrell-joint-wald}
\end{equation}
For every $S\notin\cP_0$ with $0<\Gamma(S)\le1$, it also satisfies
\begin{equation}
 \E_S\tau_\alpha(M^{\mathrm{mix}})
 \le\frac{C'_{P,\epsilon_0,d}}{\Gamma(S)}
 \left[d+\log(1/\alpha)+\ell(\Gamma(S))\right].
 \label{eq:farrell-joint-gap}
\end{equation}
Finally, for fixed full-support $S\notin\cP_0$,
\begin{equation}
 0\le n\Gamma(S)-\E_S\log M_n^{\mathrm{mix}}
 \le\log\binom{n+d-1}{d-1}+\log2.
 \label{eq:farrell-joint-regret}
\end{equation}
\end{corollary}

The word ``optimal'' in the Farrell claim refers to order and to the impossibility of a uniformly smaller rate. It does not assert a sharp leading constant or a matching lower bound at every isolated alternative. Specifically, for every fixed $\alpha\in(0,1/2)$ and every level-$\alpha$ test $\tau$ having power one against the whole complement,
\begin{equation}
 \limsup_{S\notin\cP_0:\,\Gamma(S)\downarrow0}
 \frac{\Gamma(S)\E_S\tau}{\ell(\Gamma(S))}>0.
 \label{eq:farrell-lower}
\end{equation}
To apply the general Farrell lower bound~\cite[Theorem 3.2 and Corollary 3.3]{agrawal2025stopping}, Appendix~\ref{app:farrell-lower} constructs an explicit one-parameter family approaching a common boundary null whose KL projection is that same null throughout. This verifies that the lower-bound obstruction occurs in the present TV-ball problem.
Exact implementation details for the learned-event and combined processes appear in \Cref{sec:farrell-implementation}.

The statistical roles of the two components are different. Clipped Kelly retains the full information number at each fixed alternative. The event-based component supplies uniform scale adaptation near the boundary. 

\section{Exact $O(\log d)$ dynamic implementation}
\label{sec:dynamic}

This section gives exact implementations of clipped Kelly, the Farrell-optimal processes, and their mixture.

\subsection{Clipped-Kelly implementation}
\label{sec:dynamic-ck}

We first describe a static implementation that costs $O(d\log d)$ per
oracle evaluation, and then show how to maintain the same every-round
factor in $O(\log d)$ worst-case time per observation. Throughout this
section, complexity counts unit-cost exact arithmetic operations and
comparisons, and memory is measured in numeric storage words.

Write $\rho_i=U_i/P_i$ and sort
$\rho_{(1)}\leq\cdots\leq\rho_{(d)}$. Carry $P_i,U_i$ along with the
permutation and define the prefix and suffix sums
\[
    P_j^{\mathrm{pre}}=\sum_{r=1}^j P_{(r)},
    \qquad
    U_j^{\mathrm{pre}}=\sum_{r=1}^j U_{(r)},
    \qquad
    P_j^{\mathrm{suf}}=\sum_{r=j}^d P_{(r)},
    \qquad
    U_j^{\mathrm{suf}}=\sum_{r=j}^d U_{(r)}.
\]

\begin{proposition}[Sorting-and-scanning oracle]\label{prop:static-oracle}
For $U\in\Delta_d^\circ\setminus\cP_0$, use the endpoint conventions
$\rho_{(0)}=0$ and $\rho_{(d+1)}=+\infty$. The lower threshold is
\[
    a(U)=\frac{U_j^{\mathrm{pre}}}
                  {P_j^{\mathrm{pre}}-\epsilon_0}
\]
for an index $j\in\{1,\ldots,d\}$ satisfying
\[
    P_j^{\mathrm{pre}}>\epsilon_0,
    \qquad
    \rho_{(j)}\leq a(U)\leq\rho_{(j+1)}.
\]
The upper threshold is
\[
    b(U)=\frac{U_j^{\mathrm{suf}}}
                  {P_j^{\mathrm{suf}}+\epsilon_0}
\]
for an index $j\in\{1,\ldots,d\}$ satisfying
\[
    \rho_{(j-1)}\leq b(U)\leq\rho_{(j)}.
\]
Coordinates tied at a threshold may be placed on either side.
Thus, a static oracle evaluation costs $O(d\log d)$ time and $O(d)$ memory.
\end{proposition}

Naively recomputing the oracle in \eqref{eq:plugin} sorts $d$ ratios after
every observation. We now give an exact dynamic implementation of the
same fresh factor. It maintains an implicit representation of the factor:
its value at a specified category can be evaluated in constant time,
whereas explicitly listing all $d$ values takes $O(d)$ time. The
complexity bounds are arithmetic bounds; the bit cost additionally depends
on the numerical input representations and the growing counts.
This section uses the standing assumption $\epsilon_0>0$; when
$\epsilon_0=0$, no clipping thresholds are needed and the ordinary
likelihood-ratio update is simpler.

After $t$ observations define
\[
    c_i(t)=N_i(t)+\alpha_i,
    \qquad
    C_t=t+\alpha_+,
    \qquad
    z_i(t)=\frac{c_i(t)}{P_i}.
\]
Then
\[
    \frac{\widehat U_{t,i}}{P_i}=\frac{z_i(t)}{C_t},
    \qquad
    \sum_i P_i z_i(t)=C_t.
\]
Whenever $\widehat U_t\notin\cP_0$, let $a_t,b_t$ be its oracle thresholds
and set
\[
    A_t=C_ta_t,
    \qquad
    B_t=C_tb_t.
\]
For every $t\geq0$ and $x\geq0$, define
\begin{align}
    F_t(x)
    &:=\sum_i P_i(x-z_i(t))_+-\epsilon_0x,
    \label{eq:dynamic-lower}\\
    G_t(x)
    &:=\sum_i P_i(z_i(t)-x)_+-\epsilon_0x.
    \label{eq:dynamic-upper}
\end{align}

\begin{lemma}[Rescaled threshold equations]\label{lem:rescaled}
Assume $\widehat U_t\notin\cP_0$. Then $A_t$ and $B_t$ are the unique
positive solutions of
\[
    F_t(A)=0,
    \qquad
    G_t(B)=0,
\]
respectively. Moreover,
\begin{equation}\label{eq:dynamic-factor}
    E_{\widehat U_t}(i)
    =\frac{\clip(z_i(t),A_t,B_t)}{C_t}.
\end{equation}
\end{lemma}

For any predictor, the count recursion implies that, if $X_{t+1}=k$, then
\begin{equation}\label{eq:one-key}
    z_k(t+1)=z_k(t)+\frac1{P_k},
    \qquad
    z_i(t+1)=z_i(t)\quad(i\neq k).
\end{equation}
For a nonnegative query value $v$, define
\[
    P_t^<(v)=\sum_{z_i(t)<v}P_i,
    \qquad
    S_t^<(v)=\sum_{z_i(t)<v}P_i z_i(t)
             =\sum_{z_i(t)<v}c_i(t).
\]
Then
\begin{equation}\label{eq:F-prefix}
    F_t(v)=v\{P_t^<(v)-\epsilon_0\}-S_t^<(v).
\end{equation}
Similarly, let
\[
    P_t^>(v)=\sum_{z_i(t)>v}P_i,
    \qquad
    S_t^>(v)=\sum_{z_i(t)>v}P_i z_i(t),
\]
so
\begin{equation}\label{eq:G-suffix}
    G_t(v)=S_t^>(v)-v\{P_t^>(v)+\epsilon_0\}.
\end{equation}
Write $P_t^{\geq}(v)$ and $S_t^{\geq}(v)$ for the corresponding sums over
$\{i:z_i(t)\geq v\}$.

For an outside-null predictor, $F_t$ is convex, is negative for sufficiently
small positive arguments, and has a unique positive root; after that
crossing it remains positive. Thus $F_t(v)\geq0$ is a monotone predicate
over the ordered stored keys. The function $G_t$ is strictly decreasing:
its right and left slopes are $-P_t^>(x)-\epsilon_0$ and
$-P_t^{\geq}(x)-\epsilon_0$, respectively. Hence $G_t(v)\leq0$ is also
a monotone predicate.

To locate the crossings relative to the stored keys, write
\[
    z_{\min}(t)=\min_i z_i(t),
    \qquad
    z_{\max}(t)=\max_i z_i(t).
\]
Whenever $\widehat U_t\notin\cP_0$,
\[
    F_t(C_t)=G_t(C_t)
    =C_t\{\TV(\widehat U_t,P)-\epsilon_0\}>0.
\]
Also,
\[
    F_t(z_{\min}(t))=-\epsilon_0z_{\min}(t)<0,
    \qquad
    G_t(z_{\max}(t))=-\epsilon_0z_{\max}(t)<0.
\]
The root properties therefore imply
\[
    z_{\min}(t)<A_t<C_t<B_t<z_{\max}(t).
\]
In particular, the smallest satisfying stored key in each query below
exists and has a predecessor.

Maintain the distinct $z_i(t)$ as keys of a balanced binary search tree
with worst-case height $O(\log d)$, such as an AVL or red--black tree,
grouping ties in the same node. A node at key $v$ stores its group mass,
count sum, and cardinality,
\[
    p_v=\sum_{z_i(t)=v}P_i,
    \qquad
    s_v=\sum_{z_i(t)=v}c_i(t),
    \qquad
    m_v=\#\{i:z_i(t)=v\}.
\]
Separate category arrays store $P_i$, $c_i(t)$, and $z_i(t)$, so a node
does not need a list of category indices. Once the node at $z_k(t)$ is
located, removing category $k$ subtracts $P_k$ and the old $c_k(t)$ from
its group totals and decrements its cardinality, all in constant time.
The node is deleted if its cardinality becomes zero. Insertion into an
existing group adds the updated contributions and increments its
cardinality; otherwise a new node is inserted.

Each node also stores the aggregate masses of its entire subtree
$\mathcal T$,
\[
    P(\mathcal T)=\sum_{i\in\mathcal T}P_i,
    \qquad
    S(\mathcal T)=\sum_{i\in\mathcal T}P_i z_i(t)
                 =\sum_{i\in\mathcal T}c_i(t),
\]
where $i\in\mathcal T$ means that category $i$ belongs to a node in
$\mathcal T$. Each aggregate is the sum of the current group total and
the two child-subtree totals. These aggregates are repaired in constant
time per visited node during rotations, insertion, and deletion.

\begin{lemma}[Exact lower query]\label{lem:lower-query}
Assume $\widehat U_t\notin\cP_0$. Let $v$ be the smallest stored key
satisfying $F_t(v)\geq0$. If $F_t(v)=0$, then $A_t=v$. Otherwise,
\[
    A_t=\frac{S_t^<(v)}{P_t^<(v)-\epsilon_0}.
\]
The key $v$ and $A_t$ are found by one augmented-tree descent in
$O(\log d)$ time.
\end{lemma}

\begin{lemma}[Exact upper query]\label{lem:upper-query}
Assume $\widehat U_t\notin\cP_0$. Let $v$ be the smallest stored key
satisfying $G_t(v)\leq0$. If $G_t(v)=0$, then $B_t=v$. Otherwise,
\[
    B_t=\frac{S_t^{\geq}(v)}{P_t^{\geq}(v)+\epsilon_0},
\]
where the suffix includes the node at $v$. The key $v$ and $B_t$ are
found by one augmented-tree descent in $O(\log d)$ time.
\end{lemma}

Null-ball membership is also logarithmic. Since both $\widehat U_t$ and
$P$ sum to one, $\TV(\widehat U_t,P)=\sum_i(P_i-\widehat U_{t,i})_+$.
Substituting $\widehat U_{t,i}=P_i z_i(t)/C_t$ gives
\begin{equation}\label{eq:dynamic-tv}
    \TV(\widehat U_t,P)
    =\frac{C_tP_t^<(C_t)-S_t^<(C_t)}{C_t}.
\end{equation}
Thus a single prefix query decides membership by comparing
$C_tP_t^<(C_t)-S_t^<(C_t)$ with $\epsilon_0 C_t$.

\paragraph{Initialization.}
Set $c_i(0)=\alpha_i$, $C_0=\alpha_+$, and $z_i(0)=\alpha_i/P_i$.
Form and sort the initial keys, group ties, and build the augmented
balanced tree. Then evaluate \eqref{eq:dynamic-tv} for $\widehat U_0$.
If $\widehat U_0\notin\cP_0$, perform the lower and upper queries to
initialize $A_0,B_0$; otherwise record that the first factor is
identically one. These initial queries cost $O(\log d)$ in addition
to the $O(d\log d)$ construction cost.

\paragraph{Online update.}
At the start of round $t+1$, the tree and category arrays represent
$\widehat U_t$. If $\widehat U_t\in\cP_0$, its factor is identically
one; otherwise the stored thresholds give $E_{\widehat U_t}$ through
\eqref{eq:dynamic-factor}. After observing $X_{t+1}=k$, first evaluate
$E_{\widehat U_t}(k)$ using the current state and multiply the wealth by
this factor. Then remove category $k$'s old contribution from its tie
group, and update
\[
    c_k(t+1)=c_k(t)+1,
    \qquad
    z_k(t+1)=z_k(t)+\frac1{P_k},
    \qquad
    C_{t+1}=C_t+1.
\]
All other category entries remain unchanged. Insert the updated
contribution at the new key and evaluate \eqref{eq:dynamic-tv} for
$\widehat U_{t+1}$. If this predictor is outside the null ball, perform
the two threshold queries for $A_{t+1},B_{t+1}$; otherwise record that
the next factor is identically one. Thus every betting rule is chosen
from past data, exactly as in \eqref{eq:plugin}.

\begin{theorem}[Exact dynamic every-round clipped Kelly]\label{thm:dynamic}
In the unit-cost exact-arithmetic model, after $O(d\log d)$
initialization, an implicit representation of the every-round factor
$E_{\widehat U_t}$ can be maintained using $O(d)$ numeric storage words
and $O(\log d)$ worst-case time per observation. Given this
representation, evaluating the factor at a specified category takes
$O(1)$ time. Through horizon $n$ the total cost is
\[
    O(d\log d+n\log d).
\]
Let $M_n^{\mathrm{dyn}}$ denote the wealth produced by this implementation,
with $M_0^{\mathrm{dyn}}=1$. For every observation sequence and every
$n\geq0$,
\[
    M_n^{\mathrm{dyn}}=M_n^{\mathrm{CK}}.
\]
\end{theorem}

\begin{table}[ht]
    \centering
    \caption{Computational complexity of the dynamic
    algorithm, measured in unit-cost arithmetic operations}
    \label{tab:dynamic-kelly-complexity}
    \begin{tabular}{@{}ll@{}}
        \toprule
        \textbf{Operation} & \textbf{Cost} \\
        \midrule
        Initialize and build the tree
            & $O(d\log d)$ \\
             \midrule
        Evaluate the current factor at the observed category
            & $O(1)$ \\
        Remove or insert a category contribution
            & $O(\log d)$ \\
        Check null-ball membership
            & $O(\log d)$ \\
        Compute either threshold
            & $O(\log d)$ \\
        \addlinespace
        Complete one observation update
            & {\boldmath $O(\log d)$}\\
        \bottomrule
    \end{tabular}
\end{table}

\subsection{Exact implementation of the Farrell-optimal and the combined process}
\label{sec:farrell-implementation}

The learned-event process involves infinitely many accounts, and each account bets on an event selected from a pilot sample. Nevertheless, its exact implementation requires only one active account, accumulated frozen capital, and an augmented tree of category counts. We first describe the capital bookkeeping, then show how to select and evaluate the learned events without scanning all categories.
\[
K_j+w_jW_n^{(j)}+\frac{1}{j+2},
\qquad
K_j=\sum_{r<j}w_rW_{2m_r}^{(r)}.
\]
Thus the account state requires only constant storage. At the end of block $j$, update
\[
K_{j+1}=K_j+w_jW_{2m_j}^{(j)}
\]
and initialize the next account at one.

A direct implementation maintains the category counts and scans them at each dyadic boundary to select the next event. This uses $O(d)$ memory and $O(n+d\log(n+1))$ total work through time $n$, in addition to $O(d)$ initialization. Its boundary updates cost $O(d)$. The following construction reduces these updates to $O(\log d)$.

Selecting the event by a suffix query. Maintain the raw-count ratios
\[
v_i(n)=\frac{N_i(n)}{P_i}
\]
in an augmented balanced search tree. As in Section 4.1, group equal keys and store group cardinalities together with group and subtree sums of $P_i$ and $N_i(n)$. Only the observed category changes its key, by $1/P_i$, so one removal and insertion update the tree in $O(\log d)$ worst-case time. Initially all keys are zero, allowing $O(d)$ initialization.

At a dyadic time $m=m_j$, the selected event is
\[
\widehat A_j={i(m)>m}.
\]
A single strict suffix query therefore returns
\[
p_j=\sum_{v_i(m)>m}P_i,
\qquad
s_j=\sum_{v_i(m)>m}N_i(m).
\]
These two sums determine the next bet:
\[
q_j=p_j+\epsilon_0,
\qquad
\widehat\Delta_j=
\left(\frac{s_j}{m}-p_j-\epsilon_0\right)_+,
\qquad
\lambda_j=4\widehat\Delta_j.
\]
The query does not enumerate the members of $\widehat A_j$.

Evaluating membership lazily. During the next block, the event $\widehat A_j$ must remain fixed even as the counts change. We recover its membership one observed category at a time. Maintain arrays $\mathrm{stamp}[i]$ and $\mathrm{bit}[i]$, initially setting every stamp to (-1).
\[
\mathbf1{N_i(t-1)>m_jP_i},
\qquad
\mathrm{stamp}[i]=j.
\]
At this first appearance of (i) in the block, its count has not changed since the pilot ended. Hence $N_i(t-1)=N_i(m_j)$, and the stored bit equals $\mathbf1{i\in\widehat A_j}$. Subsequent appearances reuse this bit. Incrementing the block identifier invalidates the old cache without clearing either array.

The current account is then multiplied by
\[
\exp\left\{
\lambda_j\bigl(\mathrm{bit}[i]-q_j\bigr)
-\frac{\lambda_j^2}{8}
\right\},
\]
after which the count and its tree entry are updated. Membership evaluation and account bookkeeping cost $O(1)$; the tree update costs $O(\log d)$.

Boundary updates and the combined process. The first observation completes the initial pilot and selects account zero. At each subsequent dyadic endpoint, first complete the old account’s final bet and update the counts. Then freeze its capital and select the next account from the updated tree. This order ensures that each account uses exactly its prescribed pilot and betting block.

To maintain $M^{\mathrm{mix}}$, also retain the tree of smoothed keys used by clipped Kelly. Evaluate both components using the pre-observation counts, increment the shared counts once, and update both trees. At a dyadic endpoint, select the next learned event after these updates. The two component wealths determine $M_n^{\mathrm{mix}}=(M_n^{CK}+M_n^F)/2$ in constant time.

\begin{proposition}[Logarithmic-time implementation]
\label{prop:farrell-computation}
The process $M^{\mathrm{F}}$ can be maintained with $O(d)$ initialization and memory and $O(\log d)$ worst-case time per observation. The combined process in \eqref{eq:farrell-combined} can therefore be maintained with $O(d\log d)$ initialization, $O(d)$ memory, and $O(\log d)$ worst-case time per observation.
\end{proposition}

The proof in Appendix~\ref{app:farrell-computation} verifies the maintained state and the update order.

\section{Experiment}
\label{sec:expt}

We compare the proposed plug-in clipped-Kelly process with the 
least-favorable-distribution (LFD) e-process in a four-category multinomial
model.  We take
\[
P=(0.25,0.25,0.25,0.25),\qquad
\epsilon_0=\epsilon_1=0.05,
\]
with null class \(B_{\mathrm{TV}}(P,0.05)\) and, for defining the LFD
benchmark, alternative class \(B_{\mathrm{TV}}(Q,0.05)\), where
\(Q=(0.37,0.13,0.25,0.25)\).  The corresponding least-favorable pair is
\(R^\star=(0.30,0.20,0.25,0.25)\) and
\(Q^\star=(0.32,0.18,0.25,0.25)\), giving
\(E_{\mathrm{LFD}}=Q^\star/R^\star=(1.0667,0.9,1,1)\).
We generate i.i.d.\ observations from
\[
S=(0.37,0.08,0.30,0.25),
\]
which lies on the boundary of the alternative TV ball since
\(D_{\mathrm{TV}}(S,Q)=0.05\).  At each round, our method forms the
Laplace-smoothed empirical distribution
\(\widehat S_{t-1,i}=(N_i(t-1)+1)/(t+3)\) and uses the corresponding
pointwise clipped-Kelly oracle \(E_{\widehat S_{t-1}}\).  We compare this
with the fixed LFD process and with an infeasible oracle that knows \(S\).
The alternative center \(Q\) is supplied only to the LFD benchmark. The adaptive procedure uses the null center \(P\), the robustness radius, and past observations; it does not receive \(Q\) or the data-generating distribution \(S\). The oracle curve represents the performance available if \(S\) were known.

\begin{figure}[!ht]
    \centering
    \includegraphics[width=0.6\textwidth]{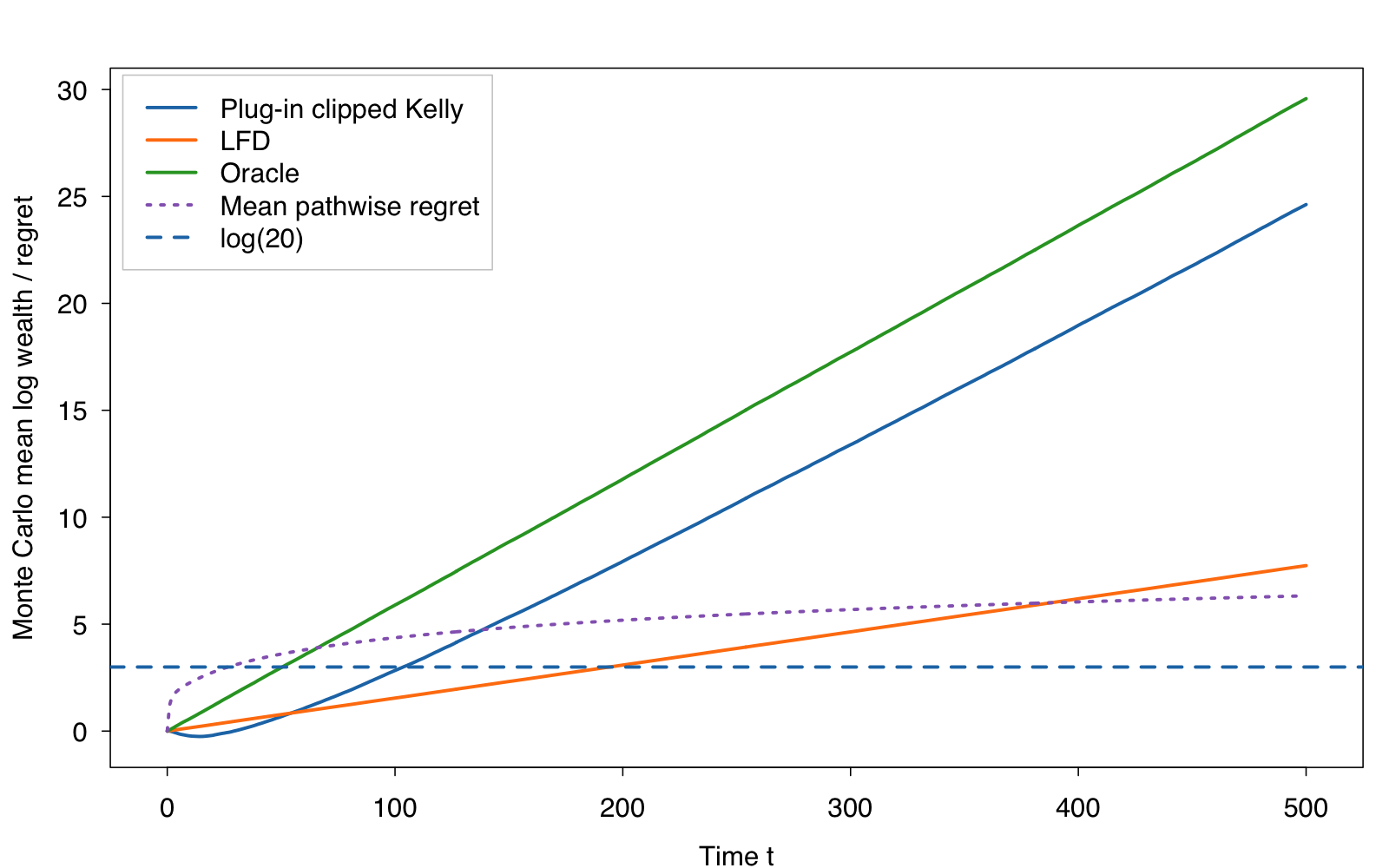}
    \caption{
    Mean cumulative log wealth in the four-category TV-neighborhood
    experiment with \(\epsilon_0=\epsilon_1=0.05\).
    The curves show the proposed predictable plug-in clipped-Kelly
    process, the fixed least-favorable-distribution (LFD) process, and
    the infeasible oracle that knows the true alternative \(S\).
    Results are averaged over \(10{,}000\) independent Monte Carlo
    replications.  Although the plug-in process
    incurs an initial learning cost, its long-run growth rate is substantially larger
    than that of the LFD process.
    }
    \label{fig:equal-eps-growth}
\end{figure}

Figure~\ref{fig:equal-eps-growth} shows the mean cumulative log wealth over
\(10{,}000\) independent simulations up to time \(500\); the dashed line is
the level-\(0.05\) threshold \(\log 20\) and the dotted line is the average pathwise regret, $\log W_t^\star-\log M_t^{\mathrm{CK}}$, of our plug-in clipped Kelly e-process. The LFD growth rate is \(0.01545\),
whereas the oracle rate is \(\Gamma(S)=0.05899\), about \(3.82\) times
larger.

The plug-in method initially incurs a learning cost but approaches
the oracle slope and substantially outgrows the fixed LFD process, with mean
stopping times \(98.1\), \(195.3\), and \(52.7\) for the plug-in, LFD, and
oracle procedures, respectively.

\section{Discussion}
\label{sec:discussion}

We studied power-one sequential testing with the goal of adapting to the unknown alternative while maintaining anytime validity over a total variation ball. Our predictable clipped-Kelly e-process learns the corresponding pointwise oracle online, with a pathwise regret bound. We also establish the optimality of the stopping time for that e-process asymptotically as the type-I error level $\alpha$ converges to zero.
We construct a complementary e-process that provides adaptation to alternatives approaching the null boundary. Averaging the two processes yields a single e-process with power one against the entire complement of the null. Its expected stopping time is simultaneously optimal in
both regimes. 

Exact dynamic implementation accompanies the statistical guarantees. After \(O(d\log d)\) initialization, the combined process requires \(O(\log d)\) worst-case time per observation. The existing LFD process has the same initialization cost and takes constant time per observation.  Thus, adaptation to the unknown alternative is achieved with only logarithmic per-observation overhead, while allowing faster evidence accumulation against the null.

The TV model serves as a concrete robust model, for which we establish statistical and computational guarantees. The statistical arguments can be extended to other convex robust null classes
on finite alphabets, including Huber's $\epsilon$-contamination
neighborhoods and general f-divergence neighborhoods. The explicit formulas and computational
guarantees, however, depend on the structure of the chosen null class.

\paragraph{Acknowledgments.} AI tools assisted with the design of computationally efficient algorithms, obtaining the regret bound, the implementation of simulation experiments, and the writing. The authors take full responsibility for all content of the paper.

\bibliographystyle{apalike}
\bibliography{ref}

\appendix

\section{Omitted Proofs}
\label{app:proofs}
\subsection{Omitted Proofs from Section 2}

\begin{proof}[Proof of Lemma \ref{lem:oracle-roots}]
For the lower threshold, define
\[
    F_-(z)=\sum_i(zP_i-U_i)_+-z\epsilon_0.
\]
It is continuous and convex.  Since $U_i/P_i>0$ for every $i$, if
$0<z<\min_i U_i/P_i$ then the positive-part sum vanishes and
$F_-(z)=-z\epsilon_0<0$.  On the other hand,
\[
    F_-(1)=\sum_i(P_i-U_i)_+-\epsilon_0
    =\TV(U,P)-\epsilon_0>0.
\]
Thus a root exists in $(0,1)$.  At any root $z$, writing $L_z=\{i:U_i/P_i<z\}$,
\[
    z\{P(L_z)-\epsilon_0\}=U(L_z)>0,
\]
so the left slope is positive.  Convexity rules out a second positive
crossing.  For the upper threshold,
\[
    F_+(z)=\sum_i(U_i-zP_i)_+-z\epsilon_0
\]
is continuous, $F_+(1)=\TV(U,P)-\epsilon_0>0$, and $F_+(z)=-z\epsilon_0<0$ once $z>\max_i U_i/P_i$.  Moreover, its right and left derivatives are, respectively,
\[
    -\sum_{i:U_i/P_i>z}P_i-\epsilon_0,
    \qquad
    -\sum_{i:U_i/P_i\geq z}P_i-\epsilon_0,
\]
which are strictly negative.  Hence $F_+$ is strictly decreasing and has a unique root in $(1,\infty)$.
\end{proof}

\begin{proof}[Proof of \Cref{thm:oracle}]
If $U\in\cP_0$, Jensen's inequality gives $\sum_i U_i\log E_i\leq\log\sum_i U_i E_i\leq0$ for every $E\in\cE_0$, so $E_U\equiv1$ is optimal and $R_U=U$.  Assume henceforth that $U\notin\cP_0$.

Let
\[
    L=\{i:U_i/P_i<a(U)\},\qquad
    H=\{i:U_i/P_i>b(U)\}.
\]
On $L$, $R_U(i)=U_i/a(U)<P_i$; on the middle set $R_U(i)=P_i$; and on
$H$, $R_U(i)=U_i/b(U)>P_i$.  By the threshold equations,
\[
    \sum_{i\in L}(P_i-R_U(i))
    =\frac1{a(U)}\sum_{i\in L}(a(U)P_i-U_i)
    =\epsilon_0,
\]
and similarly
\[
    \sum_{i\in H}(R_U(i)-P_i)=\epsilon_0.
\]
Thus $R_U$ is a probability vector and $\TV(R_U,P)=\epsilon_0$.

Next,
\begin{align*}
    1-P\cdot E_U
    &=\sum_{i\in L}(U_i-a(U)P_i)
      +\sum_{i\in H}(U_i-b(U)P_i)\\
    &=-a(U)\epsilon_0+b(U)\epsilon_0,
\end{align*}
so
\[
    P\cdot E_U+\epsilon_0\{b(U)-a(U)\}=1.
\]
Therefore, for every $R\in\cP_0$,
\[
    \E_R E_U\leq1,
\]
so $E_U\in\cE_0$.

It follows from part (a) and Theorem 4.1 of \cite{larsson2025numeraire} that $R_U$ is a reverse projection. Also, one can check that
for arbitrary $R\in\cP_0$,
\begin{align*}
    \KL(U\|R)-\KL(U\|R_U)
    &=-\sum_i U_i\log\frac{R_i}{R_U(i)}\\
    &\geq-\log\sum_i U_i\frac{R_i}{R_U(i)}\\
    &=-\log\sum_i R_i E_U(i)\geq0.
\end{align*}
Hence $R_U$ is a reverse projection.  Because $U_i>0$ for all $i$, every $R$ with finite $\KL(U\|R)$ has positive coordinates, and $R\mapsto\KL(U\|R)$ is strictly convex there; the reverse projection is therefore unique.
Finally, for any $E\in\cE_0$,
\begin{align*}
    \sum_i U_i\log E_i
    &=\KL(U\|R_U)+\sum_i U_i\log\frac{E_i}{E_U(i)}\\
    &\leq\KL(U\|R_U)+\log\sum_i R_U(i)E_i\\
    &\leq\KL(U\|R_U).
\end{align*}
Equality holds for $E=E_U$, and strictness in Jensen together with full
support of $U$ gives uniqueness.
\end{proof}

For $m,\delta>0$, define
\[
    \cK_{m,\delta}
    :=\{U\in\Delta_d:\min_i U_i\geq m,
      \ \TV(U,P)\geq\epsilon_0+\delta\}.
\]

\begin{lemma}[Continuity and uniform boundedness]\label{lem:continuity}
On every nonempty $\cK_{m,\delta}$, the maps
$U\mapsto a(U)$, $U\mapsto b(U)$, $U\mapsto E_U$, and
$U\mapsto\Gamma(U)$ are continuous.  Moreover,
\[
    \frac{m}{\max_i P_i}
    \leq E_U(i)\leq\frac1{\min_i P_i}
\]
uniformly over $U\in\cK_{m,\delta}$ and $i\leq d$.
\end{lemma}

\begin{proof}
For the lower root define
\[
    F(U,z)=\sum_i(zP_i-U_i)_+-z\epsilon_0.
\]
This function is jointly continuous.  For $U\in\cK_{m,\delta}$, the unique root satisfies
\[
    \frac{m}{\max_i P_i}\leq a(U)<1,
\]
and $F(U,1)=\TV(U,P)-\epsilon_0\geq\delta$, so no sequence of lower roots can converge to one.  If $U_n\to U$, compactness gives a convergent subsequence $a(U_{n_k})\to a_\infty$; joint continuity yields $F(U,a_\infty)=0$, and uniqueness forces $a_\infty=a(U)$.  Hence the full sequence converges.  The same argument applies to $b(U)$, whose roots lie in the common compact interval $[1,1/\min_i P_i]$ and are separated from one by the same margin condition.

Clipping is continuous, so $U\mapsto E_U$ is continuous.  Since clipping keeps every coordinate between the minimum and maximum of $U_i/P_i$,
\[
    \frac{m}{\max_i P_i}\leq E_U(i)\leq\frac1{\min_i P_i}.
\]
The lower bound is positive uniformly on $\cK_{m,\delta}$, and therefore $\Gamma(U)=\sum_i U_i\log E_U(i)$ is continuous as well.
\end{proof}

\begin{proof}[Proof of Lemma \ref{lem:oracle-mismatch}]
The lower bound follows because $E_U\in\cE_0$ and $E_S$ is pointwise optimal at $S$.  For the upper bound, recall that $E_U=U/R_U$ for some $R_U\in\cP_0$.  Hence
\[
    \KL(S\|R_U)
    =\KL(S\|U)+h(S,U).
\]
Since $R_S$ minimizes $R\mapsto\KL(S\|R)$ over $\cP_0$,
\[
    \Gamma(S)=\KL(S\|R_S)\leq\KL(S\|R_U),
\]
which rearranges to the upper bound in \eqref{eq:oracle-mismatch}.
\end{proof}

\begin{proof}[Proof of \Cref{thm:adaptive-valid}]
Fix $R\in\cP_0$ and suppose that $X_1,X_2,\ldots$ are i.i.d.\
with distribution $R$.
For every $t\geq1$, the predictor $\widehat U_{t-1}$ is
$\mathcal F_{t-1}$-measurable and has strictly positive coordinates.
By Theorem~\ref{thm:oracle},
$E_{\widehat U_{t-1}}\in\cE_0$, and hence
\[
    \sum_{i=1}^d R_i E_{\widehat U_{t-1}}(i)\leq1.
\]
Since $X_t$ is independent of $\mathcal F_{t-1}$ with
distribution $R$, we obtain
\[
    \begin{aligned}
    \mathbb E_R\!\left[
        M_t^{\mathrm{CK}}\mid\mathcal F_{t-1}
    \right]
    &=
    M_{t-1}^{\mathrm{CK}}
    \mathbb E_R\!\left[
        E_{\widehat U_{t-1}}(X_t)\mid\mathcal F_{t-1}
    \right]\\
    &=
    M_{t-1}^{\mathrm{CK}}
    \sum_{i=1}^d R_i E_{\widehat U_{t-1}}(i)\\
    &\leq M_{t-1}^{\mathrm{CK}}.
    \end{aligned}
\]
The process is adapted and nonnegative, with
$M_0^{\mathrm{CK}}=1$.
Taking expectations and iterating gives
$\mathbb E_R M_t^{\mathrm{CK}}\leq1$ for every $t$,
which also establishes integrability.
Thus $(M_t^{\mathrm{CK}})_{t\geq0}$ is a nonnegative test
supermartingale under every i.i.d.\ law $R^\infty$
with $R\in\cP_0$.
\end{proof}
Now we state and prove the following auxiliary lemma, which we will invoke later.
\begin{lemma}[Logarithmic envelope]\label{lem:envelope}
There is a constant $C<\infty$, depending only on $P$ and $\alpha$, such
that every smoothed empirical predictor based on $T$ observations satisfies
\[
    \max_i|\log E_U(i)|\leq C+\log(T+1).
\]
\end{lemma}
\begin{proof}[Proof of Lemma \ref{lem:envelope}]
Every coordinate satisfies
$U_i\geq\alpha_{\min}/(T+\alpha_+)$ and $U_i\leq1$.  If $U\in\cP_0$,
then $E_U\equiv1$.  Otherwise, clipping keeps $E_U(i)$ between the minimum
and maximum of $U_i/P_i$.
\end{proof}

\begin{proof}[Proof of \Cref{thm:plugin-regret}]
Fix $n\geq1$ and a sequence $x^n$, and write $N_i=N_i(n)$.
We use $\log0=-\infty$, omit zero-weight terms in logarithmic
sums, and use $0^0=1$ in likelihood products.

For $\mathcal A=\cE_0$ or $\mathcal A=\cP_0$, define
\[
    J_{\mathcal A}(c)
    :=\max_{z\in\mathcal A}\sum_i c_i\log z_i,
    \qquad c\in[0,\infty)^d.
\]
These maxima exist and are finite: both classes are compact,
the objective is upper semicontinuous, and each class contains
a vector with strictly positive coordinates.
Indeed, $\cE_0$ is closed, and since $P\in\cP_0$,
every $E\in\cE_0$ satisfies
\[
    0\leq E_i\leq \frac{1}{P_i},
    \qquad i=1,\ldots,d.
\]
Thus $\cE_0$ is compact; compactness of $\cP_0$ follows from its
definition as a closed subset of the simplex.

Let $c_t=\one+N(t)$ and choose a maximizer $z_t$ at $c_t$.
Since every coordinate of $c_t$ is positive, every coordinate
of $z_t$ is positive. Evaluating the objective at the previous
maximizer gives
\[
    J_{\mathcal A}(c_t)
    \geq
    \sum_i c_{t,i}\log z_{t-1,i}
    =
    J_{\mathcal A}(c_{t-1})+\log z_{t-1}(x_t).
\]
Moreover,
\[
    J_{\mathcal A}(c+c')
    \leq J_{\mathcal A}(c)+J_{\mathcal A}(c'),
    \qquad c,c'\in[0,\infty)^d.
\]
Consequently, summing over $t$ yields
\[
    \sum_{t=1}^n\log z_{t-1}(x_t)
    \leq
    J_{\mathcal A}(\one+N)-J_{\mathcal A}(\one)
    \leq J_{\mathcal A}(N).
\]
By Corollary~\ref{cor:ftrl}, we may take
$z_t=E_{\widehat U_t}$ when $\mathcal A=\cE_0$.
Therefore,
\[
    \log M_n^{\mathrm{CK}}
    \leq J_{\cE_0}(N)
    =\log W_n^\star.
\]
This proves the nonnegativity of the regret.

Next, set $R_t=R_{\widehat U_t}$ and define
\[
    D_n(x^n):=\prod_{t=1}^n R_{t-1}(x_t),
    \qquad
    L_0(x^n):=\max_{R\in\cP_0}\prod_iR_i^{N_i}.
\]
Since $\widehat U_t=c_t/(t+d)$ has full support,
Theorem~\ref{thm:oracle} implies that $R_t$ maximizes
$\sum_i c_{t,i}\log R_i$ over $\cP_0$.
Applying the preceding telescoping inequality to $\cP_0$
therefore gives
\[
    \log D_n(x^n)
    \leq J_{\cP_0}(N)
    =\log L_0(x^n).
\]
Also define the unrestricted maximum likelihood
\[
    L_\Delta(x^n)
    :=\max_{q\in\Delta_d}\prod_iq_i^{N_i}
    =\prod_i\left(\frac{N_i}{n}\right)^{N_i}.
\]

We now establish
\[
    \log W_n^\star
    =
    \log L_\Delta(x^n)-\log L_0(x^n)
    =
    n\Gamma(\bar S_n),
\]
including when $\bar S_n$ lies on the boundary of the simplex.
Let
\[
    I:=\{i:N_i>0\},
\]
and let $R^\star$ attain $L_0(x^n)$.
Since $P\in\cP_0$ has full support, $L_0(x^n)>0$;
hence $R_i^\star>0$ for every $i\in I$.

For every $E\in\cE_0$, Jensen's inequality and null validity
imply
\[
    \sum_{i\in I}\bar S_{n,i}
       \log\frac{R_i^\star E_i}{\bar S_{n,i}}
    \leq
    \log\left(\sum_{i\in I}R_i^\star E_i\right)
    \leq0,
\]
where the inequality is immediate if $E_i=0$ for some $i\in I$.
Rearranging and maximizing over $E\in\cE_0$ gives
\[
    \log W_n^\star
    \leq
    \log L_\Delta(x^n)-\log L_0(x^n).
\]

To prove the reverse inequality, define
\[
    E_i^\star:=
    \begin{cases}
        \bar S_{n,i}/R_i^\star, & i\in I,\\
        0, & i\notin I.
    \end{cases}
\]
We claim that $E^\star\in\cE_0$.
Indeed, for every $R\in\cP_0$, convexity of $\cP_0$ ensures that
\[
    R_\lambda:=(1-\lambda)R^\star+\lambda R
    \in\cP_0,
    \qquad 0\leq\lambda\leq1.
\]
Because $R^\star$ maximizes
$\sum_{i\in I}\bar S_{n,i}\log R_i$ over $\cP_0$,
the right derivative at $\lambda=0$ satisfies
\[
    0
    \geq
    \left.
    \frac{d}{d\lambda}
    \sum_{i\in I}\bar S_{n,i}\log R_{\lambda,i}
    \right|_{\lambda=0+}
    =
    \sum_{i\in I}
        \bar S_{n,i}\frac{R_i-R_i^\star}{R_i^\star}.
\]
Since $\sum_{i\in I}\bar S_{n,i}=1$, this yields
\[
    \sum_i R_iE_i^\star
    =
    \sum_{i\in I}\frac{R_i\bar S_{n,i}}{R_i^\star}
    \leq1.
\]
Thus $E^\star\in\cE_0$, and
\[
    \sum_{i\in I}N_i\log E_i^\star
    =
    \log L_\Delta(x^n)-\log L_0(x^n).
\]
Together with the preceding upper bound, this proves
\[
    \log W_n^\star
    =
    \log L_\Delta(x^n)-\log L_0(x^n).
\]
Furthermore, maximizing the null likelihood is equivalent to
minimizing $D_{\mathrm{KL}}(\bar S_n\|R)$ over $R\in\cP_0$.
Therefore,
\[
    \log W_n^\star
    =
    nD_{\mathrm{KL}}(\bar S_n\|R^\star)
    =
    n\Gamma(\bar S_n).
\]

Finally, define
\[
    Q_{\mathrm L}(x^n)
    :=\prod_{t=1}^n\widehat U_{t-1,x_t}.
\]
The identity
$E_{\widehat U_{t-1}}=\widehat U_{t-1}/R_{t-1}$
gives
\[
    M_n^{\mathrm{CK}}
    =\frac{Q_{\mathrm L}(x^n)}{D_n(x^n)}.
\]
Using $\log D_n(x^n)\leq\log L_0(x^n)$, we obtain
\begin{align*}
    \log W_n^\star-\log M_n^{\mathrm{CK}}
    &=
    \log L_\Delta(x^n)-\log L_0(x^n)
    -\log Q_{\mathrm L}(x^n)+\log D_n(x^n)\\
    &\leq
    \log\frac{L_\Delta(x^n)}{Q_{\mathrm L}(x^n)}.
\end{align*}
Consequently,
\begin{equation}\label{eq:regret-to-predictive-kl}
    0\leq\log W_n^\star-\log M_n^{\mathrm{CK}}
    \leq\log\frac{L_\Delta(x^n)}{Q_{\mathrm L}(x^n)}.
\end{equation}

For the Laplace predictor,
\begin{equation}\label{eq:laplace-sequence-prob}
    Q_{\mathrm L}(x^n)
    =
    \frac{(d-1)!\prod_{i=1}^d N_i!}{(n+d-1)!}.
\end{equation}
Moreover,
\[
    \frac{n!}{\prod_iN_i!}
    \prod_i\left(\frac{N_i}{n}\right)^{N_i}
    \leq1,
\]
because the left-hand side is the probability, under the
multinomial law with parameter $\bar S_n$, of observing
the count vector $N$.
Therefore,
\begin{equation}\label{eq:laplace-ratio-bound}
    \frac{L_\Delta(x^n)}{Q_{\mathrm L}(x^n)}
    \leq
    \frac{(n+d-1)!}{n!(d-1)!}
    =
    \binom{n+d-1}{d-1}.
\end{equation}
Combining \eqref{eq:regret-to-predictive-kl},
\eqref{eq:laplace-ratio-bound}, and
$\log W_n^\star=n\Gamma(\bar S_n)$ proves the result.
\end{proof}
\begin{proof}[Proof of Corollary \ref{cor:plugin-expected-growth}]
Since $E_S\in\cE_0$, Theorem~\ref{thm:plugin-regret} implies, pathwise,
\[
    \sum_{t=1}^n\log E_S(X_t)-\log M_n^{\mathrm{CK}}
    \leq\log W_n^\star-\log M_n^{\mathrm{CK}}
    \leq\log\binom{n+d-1}{d-1}.
\]
Taking expectations under $S$ gives the upper bound in
\eqref{eq:plugin-regret}, since $\E_S\log E_S(X_t)=\Gamma(S)$.
The lower bound follows from conditional pointwise optimality:
$\E_S[\log E_{\widehat U_{t-1}}(X_t)\mid\cF_{t-1}]
=h(S,\widehat U_{t-1})\leq\Gamma(S)$.
The simpler bound \eqref{eq:plugin-regret-simple} follows from
$\binom{N}{k}\leq(eN/k)^k$.
Dividing \eqref{eq:plugin-regret} by $n$ proves
\eqref{eq:plugin-expected-growth}.
\end{proof}

\begin{proof}[Proof of Corollary \ref{cor:plugin-as-growth}]
Put
\[
    \Delta_t
    :=\Gamma(S)-h(S,\widehat U_{t-1})\geq0.
\]
Since $\E_S\sum_{t=1}^n\Delta_t
=n\Gamma(S)-\E_S\log M_n^{\mathrm{CK}}$,
Corollary~\ref{cor:plugin-expected-growth} gives
\[
    \E_S\sum_{t=1}^n\Delta_t
    \leq \log\binom{n+d-1}{d-1}.
\]
For $n_k=2^k$, Markov's inequality yields, for every $\varepsilon>0$,
\[
    \Pp_S\!\left(
      \sum_{t=1}^{n_k}\Delta_t>\varepsilon n_k
    \right)
    \leq
    \frac{\log\binom{2^k+d-1}{d-1}}
         {\varepsilon 2^k}.
\]
The right-hand side is summable in $k$.  Borel--Cantelli therefore gives
$n_k^{-1}\sum_{t\leq n_k}\Delta_t\to0$ almost surely; monotonicity of the
partial sums extends this to all $n$.

Now define
\[
    D_t
    :=\log E_{\widehat U_{t-1}}(X_t)
      -h(S,\widehat U_{t-1}).
\]
Then $(D_t)$ is a martingale-difference sequence.  By
Lemma~\ref{lem:envelope},
$|D_t|\leq C'(1+\log(t+1))$, so
\[
    \sum_{t=1}^\infty
    \frac{\E_S[D_t^2\mid\cF_{t-1}]}{t^2}<\infty
    \qquad\text{almost surely}.
\]
The martingale strong law and Kronecker's lemma imply
$n^{-1}\sum_{t\leq n}D_t\to0$ almost surely.  Since
\[
    \frac1n\log M_n^{\mathrm{CK}}
    =\Gamma(S)
      -\frac1n\sum_{t=1}^n\Delta_t
      +\frac1n\sum_{t=1}^nD_t,
\]
\eqref{eq:plugin-as-growth} follows.
\end{proof}

\begin{proof}[Proof of \Cref{thm:stop-lower}]
If $\E_S\tau=\infty$ there is nothing to prove.  Put $\tau_m=\tau\wedge m$ and write
\[
    Y_t=\log\frac{S(X_t)}{R_S(X_t)},
    \qquad \E_sY_t=\Gamma(S).
\]
The likelihood ratio of the stopped observation under $S$ relative to $R_S$ is $\prod_{t=1}^{\tau_m}S(X_t)/R_S(X_t)$, and hence
\[
    \KL(\Pp_S^{\tau_m}\|\Pp_{R_S}^{\tau_m})
    =\E_S\sum_{t=1}^{\tau_m}Y_t.
\]
Because $\tau_m\leq m$ and $\{\tau_m\geq t\}\in\cF_{t-1}$, independence under the i.i.d. law $S$ gives
\begin{align*}
    \E_S\sum_{t=1}^{\tau_m}Y_t
    &=\sum_{t=1}^m\E_S[\one\{\tau_m\geq t\}Y_t]\\
    &=\Gamma(S)\sum_{t=1}^m\Pp_S(\tau_m\geq t)\\
    &=\Gamma(S)\E_S\tau_m.
\end{align*}
By data processing through the indicator $\one\{\tau\leq m\}$,
\[
    \Gamma(S)\E_S\tau_m
    \geq\KL(p_m,q_m),
\]
where $p_m=\Pp_S(\tau\leq m)\to1$ and
$q_m=\Pp_{R_S}(\tau\leq m)\leq\alpha$.  Hence
$\liminf_m\KL(p_m,q_m)\geq\log(1/\alpha)$.  Monotone convergence gives
\eqref{eq:stop-lower}.
\end{proof}

\begin{proof}[Proof of \Cref{thm:stop-upper}]

Write $\Gamma=\Gamma(S)>0$ and fix $0<\eta<\Gamma/8$.  By
Lemma~\ref{lem:continuity}, there are $\delta,L>0$ such that
\begin{equation}\label{eq:good-neighborhood}
    \|U-S\|_1\leq\delta
    \quad\Longrightarrow\quad
    h(S,U)\geq\Gamma-\eta,
    \qquad
    \max_i|\log E_U(i)|\leq L.
\end{equation}
Let $\bar S_t=N(t)/t$ for $t\geq1$.  For the Laplace predictor,
\[
    \|\widehat U_t-\bar S_t\|_1
    \leq \frac{2d}{t+d}.
\]
For all sufficiently large $t$ this is at most $\delta/2$.  If $\|\bar S_t-S\|_1>\delta/2$, then some coordinate differs from its mean by more than $\delta/(2d)$.  Hoeffding's inequality and a union bound therefore give
\begin{equation}\label{eq:predictor-tail}
    \Pp_S(\|\widehat U_t-S\|_1>\delta)
    \leq 2d\exp\!\left(-\frac{t\delta^2}{2d^2}\right)
\end{equation}
for all sufficiently large $t$.

For a deterministic horizon $n$, let $T_n=\lceil K\log n\rceil$.  Summing \eqref{eq:predictor-tail} over $t\geq T_n$ and choosing $K$ sufficiently large gives constants $C<\infty$ and $q>2$ such that
\[
    \Pp_S\!\left(
      \exists t\geq T_n:\|\widehat U_t-S\|_1>\delta
    \right)
    \leq Cn^{-q}.
\]  Define
\[
    Z_t=\log E_{\widehat U_{t-1}}(X_t),
    \qquad
    m_t=\E_S[Z_t\mid\cF_{t-1}],
\]
and let $\sigma$ be the first $t>T_n$ for which the predictor
$\widehat U_{t-1}$ leaves the $\delta$-neighborhood of $S$.  The stopped
increments
\[
    (Z_t-m_t)\one\{t<\sigma\},\qquad t>T_n,
\]
form a martingale-difference sequence bounded in absolute value by $2L$.
Azuma--Hoeffding therefore yields
\begin{equation}\label{eq:post-burn-tail}
    \Pp_S\!\left(
       \sigma>n,\ 
       \sum_{t=T_n+1}^n Z_t
       <(n-T_n)(\Gamma-2\eta)
    \right)
    \leq e^{-c_\eta(n-T_n)}.
\end{equation}

The first $T_n$ increments satisfy, deterministically by
Lemma~\ref{lem:envelope},
\[
    \sum_{t=1}^{T_n}Z_t
    \geq-O(T_n\log T_n)
    =-o(n).
\]
Combining this with \eqref{eq:predictor-tail} and
\eqref{eq:post-burn-tail}, for all sufficiently large $n$,
\begin{equation}\label{eq:linear-growth-tail}
    \Pp_S\!\left(
       \log M_n^{\mathrm{CK}}<n(\Gamma-3\eta)
    \right)
    \leq Cn^{-q}+e^{-c_\eta n}.
\end{equation}

Put $B=\log(1/\alpha)$ and
\[
    n_B=\left\lceil\frac{B}{\Gamma-4\eta}\right\rceil.
\]
For all sufficiently large $B$ and every $n\geq n_B$,
$n(\Gamma-3\eta)>B$.  Therefore
\[
    \Pp_S(\tau_\alpha^{\mathrm{CK}}>n)
    \leq Cn^{-q}+e^{-c_\eta n}.
\]
Using $\E\tau=\sum_{n\geq0}\Pp(\tau>n)$ gives
\begin{align*}
    \E_S\tau_\alpha^{\mathrm{CK}}
    &\leq n_B+\sum_{n\geq n_B}\{Cn^{-q}+e^{-c_\eta n}\}\\
    &=n_B+o(B),
\end{align*}
because $q>2$ and $n_B\asymp B$.
Hence
\[
    \limsup_{\alpha\downarrow0}
    \frac{\E_S\tau_\alpha^{\mathrm{CK}}}{\log(1/\alpha)}
    \leq\frac1{\Gamma-4\eta}.
\]
Let $\eta\downarrow0$.  The matching lower bound follows from
Theorem~\ref{thm:stop-lower}, since anytime validity implies
$\Pp_{R_S}(\tau_\alpha^{\mathrm{CK}}<\infty)\leq\alpha$ and
Corollary~\ref{cor:plugin-as-growth} implies $\tau_\alpha^{\mathrm{CK}}<\infty$ $S$-almost surely.
\end{proof}

\subsection{Omitted proofs from Section 3}

\begin{proof}[Proof of Lemma~\ref{lem:farrell-geometry}]
For any $R\in\cP_0$, the triangle inequality gives $\TV(S,R)\ge\Delta(S)$. Pinsker's inequality therefore gives $\KL(S\|R)\ge2\Delta(S)^2$, proving the lower bound.

For the upper bound first suppose $S$ has full support. Let $a<1<b$ be the clipped-Kelly thresholds from \Cref{thm:oracle}, and write
\[
 L=\{i:S_i/P_i<a\},\quad H=\{i:S_i/P_i>b\},\quad
 E_S(i)=\clip(S_i/P_i,a,b),\quad R_S(i)=S_i/E_S(i).
\]
The threshold equations imply
\[
 a=\frac{S(L)}{P(L)-\epsilon_0},\qquad
 b=\frac{S(H)}{P(H)+\epsilon_0},\qquad P(L)>\epsilon_0.
\]
Every event discrepancy has absolute value at most $\epsilon_0+\Delta(S)$. Consequently,
\[
 1-a\le\frac{\Delta(S)}{P(L)-\epsilon_0}\le\frac{\Delta(S)}{\eta_0},
 \qquad
 b-1\le\frac{\Delta(S)}{P(H)+\epsilon_0}\le\frac{\Delta(S)}{\epsilon_0}.
\]
Hence $|E_S(i)-1|\le\Delta(S)/c_0$ for every $i$. Since $S_i=R_S(i)E_S(i)$ and both vectors sum to one,
\[
 \Gamma(S)=\KL(S\|R_S)
 \le\sum_i\frac{(S_i-R_S(i))^2}{R_S(i)}
 =\sum_iR_S(i)(E_S(i)-1)^2
 \le c_0^{-2}\Delta(S)^2.
\]

To extend the result to boundary alternatives, we verify continuity of $\Gamma$ on the entire simplex. Compactness of $\cP_0$ and joint lower semicontinuity of KL give lower semicontinuity. For upper semicontinuity at $S$, choose a minimizer $R^*\in\cP_0$ and let $R_\rho=(1-\rho)R^*+\rho P$ for $0<\rho<1$. This is a strictly positive null vector. If $S_k\to S$, then
\[
 \limsup_k\Gamma(S_k)\le\KL(S\|R_\rho)
 \le\Gamma(S)-\log(1-\rho).
\]
Letting $\rho\downarrow0$ proves continuity. Approximation by full-support alternatives completes the proof. 

\end{proof}

\begin{proof}[Proof of Theorem~\ref{thm:farrell-upper}]

\textbf{Validity.}
For any fixed event $A$, put $q_A=P(A)+\epsilon_0$. If $R\in\cP_0$ and $\lambda\ge0$, Hoeffding's lemma gives
\[
 \E_R\exp\{\lambda(\one\{X\in A\}-q_A)-\lambda^2/8\}\le1,
\]
since $R(A)\le q_A$. The same inequality holds conditionally for predictable $A$ and $\lambda$. Each $W^{(j)}$ is constant before its block, uses predictable null-valid factors during the block, and is constant afterwards. Thus it is a test supermartingale. Conditional Tonelli applies to their nonnegative mixture. The mixture starts at one and is integrable because its expectation is at most one.\\

\textbf{An accurate pilot gives a useful bet.}
Write $\Delta=\Delta(S)>0$ and let
\[
 G_j=\{\TV(\bar S_{m_j},S)\le\Delta/8\}.
\]
Using $\TV(\bar S_m,S)=\max_A\{\bar S_m(A)-S(A)\}$, a union bound over $2^d$ events and one-sided Hoeffding give
\begin{equation}
 \Pp_S(G_j^c)\le 2^d\exp(-m_j\Delta^2/32).
 \label{eq:farrell-pilot-failure}
\end{equation}
This is only an analysis step: the algorithm does not enumerate events. On $G_j$,
\[
 \frac78\Delta\le\widehat\Delta_j\le\frac98\Delta,
 \qquad
 S(\widehat A_j)-q_j\ge\widehat\Delta_j-\frac18\Delta.
\]
In particular, conditional on the pilot, the mean log increment in the next block satisfies
\[
 \begin{split}
 \mu_j&=4\widehat\Delta_j(S(\widehat A_j)-q_j)-2\widehat\Delta_j^2\\
 &\ge2\widehat\Delta_j(\widehat\Delta_j-\Delta/4)
 \ge\frac{35}{32}\Delta^2\ge\Delta^2,
 \qquad \lambda_j\le\frac92\Delta.
 \end{split}
\]
The log increment has range length $\lambda_j$. Independence of the next block and the pilot, followed by Hoeffding's inequality, gives on $G_j$
\[
 \Pp_S\left(\left.\log W_{2m_j}^{(j)}<\tfrac12m_j\Delta^2
                  \,\right|\cF_{m_j}\right)
 \le\exp\left(-\frac{m_j\Delta^4}{2\lambda_j^2}\right)
 \le\exp(-m_j\Delta^2/41).
\]
Combining with \eqref{eq:farrell-pilot-failure},
\begin{equation}
 \Pp_S\left(\log W_{2m_j}^{(j)}<\tfrac12m_j\Delta^2\right)
 \le(2^d+1)\exp(-m_j\Delta^2/41).
 \label{eq:farrell-account-failure}
\end{equation}

\textbf{Control the entire stopping-time tail.}
Let $\kappa_j=\log((j+1)(j+2)/\alpha)$ and let $j_*$ be the first nonnegative integer with
\begin{equation}
 2^{j_*}\Delta^2\ge
 \max\{2\kappa_{j_*},\ 41[\log(2^d+1)+2]\}.
 \label{eq:farrell-epoch}
\end{equation}
Such an integer exists. If $m_j\Delta^2/2\ge \kappa_j$, then $M_{2m_j}^{\mathrm{F}}\ge w_jW_{2m_j}^{(j)}$ and \eqref{eq:farrell-account-failure} imply
\[
 \Pp_S(\tau_\alpha(M^{\mathrm{F}})>2m_j)
 \le(2^d+1)\exp(-m_j\Delta^2/41).
\]
Both conditions in \eqref{eq:farrell-epoch} persist at subsequent dyadic times: the left side doubles, whereas $\kappa_{j+1}-\kappa_j=\log((j+3)/(j+1))\le\log3$.

Put $m=2^{j_*}$, $B_d=2^d+1$, and $t=m\Delta^2/41\ge\log B_d+2$. Summing the hitting-time tail over dyadic intervals yields
\[
 \begin{split}
 \E_S\tau_\alpha(M^{\mathrm{F}})
 &\le2m+2m\sum_{r\ge0}2^r
             \Pp_S(\tau_\alpha(M^{\mathrm{F}})>2m2^r)\\
 &\le2m+2mB_d\sum_{r\ge0}2^re^{-t2^r}
 \le2m+4mB_de^{-t}<3m.
 \end{split}
\]
For the penultimate inequality, use $2^r\ge r+1$ and $t\ge2$ to bound the series by $e^{-t}/(1-2e^{-t})\le2e^{-t}$.

\textbf{Solve for the epoch and convert to information.}
Set $u=\log(1/\Delta^2)\ge0$ and
\[
 H=d+\log(1/\alpha)+\log(e+u)\ge3.
\]
Choose $2^j$ to be the first dyadic integer at least $256H/\Delta^2$. Then $2^j<512H/\Delta^2$ and
\[
 j+2\le11+\frac{u+\log H}{\log2}\le8H(e+u),
 \qquad \log(j+2)\le3H.
\]
It follows that $2\kappa_j\le14H$, while $41[\log(2^d+1)+2]\le100d\le100H$. Thus this $j$ satisfies \eqref{eq:farrell-epoch}, and
\[
 \E_S\tau_\alpha(M^{\mathrm{F}})\le3\,2^{j_*}
 <1536H/\Delta^2.
\]
This is \eqref{eq:farrell-delta-bound}. By Lemma~\ref{lem:farrell-geometry},
\[
 \Delta^{-2}\le c_0^{-2}\Gamma(S)^{-1},\qquad
 \log(1/\Delta^2)\le\log(1/\Gamma(S))+\log(c_0^{-2}).
\]
For $\Gamma(S)\le1$, the last fixed additive constant changes $\ell(\Gamma(S))$ by at most a constant depending on $c_0$. Absorbing it proves \eqref{eq:farrell-gamma-bound}.

Finally, \eqref{eq:farrell-account-failure} is summable in $j$. By Borel--Cantelli, all sufficiently large frozen accounts have $W_{2m_j}^{(j)}\ge\exp(m_j\Delta^2/2)$ almost surely. Their weighted frozen contributions tend to infinity and are retained forever. The other contributions are nonnegative, so $M_n^{\mathrm{F}}\to\infty$ almost surely. 

\end{proof}

\begin{proof}[Proof of Corollary~\ref{cor:farrell-joint}]
Mixtures preserve the test-supermartingale property. Also,
\[
 \tau_\alpha(M^{\mathrm{mix}})\le
 \tau_{\alpha/2}(M^{\mathrm{CK}})\wedge\tau_{\alpha/2}(M^{\mathrm{F}})
 \quad\text{pathwise}.
\]
Theorem~\ref{thm:farrell-upper} proves power one and \eqref{eq:farrell-joint-gap}. For a fixed full-support $S$, \Cref{thm:stop-upper} and $\log(2/\alpha)\sim\log(1/\alpha)$ give the Wald upper bound. The information lower bound applies to this valid, power-one test and gives $\E_S\tau_\alpha(M^{\mathrm{mix}})\ge\log(1/\alpha)/\Gamma(S)$, proving \eqref{eq:farrell-joint-wald}.

For regret, $M_n^{\mathrm{mix}}\ge M_n^{\mathrm{CK}}/2$ gives the upper bound in \eqref{eq:farrell-joint-regret}. To prove the lower bound, let $L_n=\prod_{t=1}^n S(X_t)/R_S(X_t)$. Validity under $R_S$ and Jensen's inequality give
\[
 \E_S\log M_n^{\mathrm{mix}}
 =n\Gamma(S)+\E_S\log(M_n^{\mathrm{mix}}/L_n)
 \le n\Gamma(S)+\log\E_{R_S}M_n^{\mathrm{mix}}
 \le n\Gamma(S).
\]
The nonnegativity statement here is for expected regret. The original nonnegative \emph{pathwise} regret need not transfer to an arbitrary mixture.

\end{proof}

\subsection{Omitted Proofs from Section 4}

\begin{proof}[Proof of Proposition \ref{prop:static-oracle}]
At the lower root the active set is a prefix, so
\[
    aP_j^{\mathrm{pre}}-U_j^{\mathrm{pre}}=a\epsilon_0.
\]
The displayed formula follows, and the interval condition is exactly
self-consistency of that prefix. The upper active set is a suffix, giving
\[
    U_j^{\mathrm{suf}}-bP_j^{\mathrm{suf}}=b\epsilon_0.
\]
Coordinates at a threshold contribute zero to its defining equation, so
including or excluding them does not change the result. After sorting,
the cumulative sums and a scan identify both thresholds in $O(d)$ time.
Sorting therefore dominates the computational cost.
\end{proof}

\begin{proof}[Proof of Lemma \ref{lem:rescaled}]
Substitute $\widehat U_{t,i}=P_i z_i(t)/C_t$ into
\eqref{eq:oracle-low}--\eqref{eq:oracle-high} and multiply by $C_t$.
This gives $F_t(A_t)=G_t(B_t)=0$, and uniqueness follows from the
corresponding oracle threshold equations.
Equation \eqref{eq:dynamic-factor} follows by positive homogeneity of
clipping.
\end{proof}

\begin{proof}[Proof of Lemma \ref{lem:lower-query}]
The crossing inequalities ensure that $v$ exists and has a predecessor.
The positive root is $v$ if equality holds; otherwise it lies strictly
between the predecessor of $v$ and $v$. On that interval the active lower
set is exactly $\{i:z_i(t)<v\}$, so
\[
    F_t(x)=x\{P_t^<(v)-\epsilon_0\}-S_t^<(v).
\]
Because the function crosses from negative to positive on this interval,
$P_t^<(v)-\epsilon_0>0$, and solving $F_t(x)=0$ gives the formula.

For the search, carry running prefix totals $(p_0,s_0)$ for all keys
preceding the current subtree, initially $(0,0)$. At a node with key $w$
and left subtree $\mathcal L$, compute
\[
    p=p_0+P(\mathcal L),
    \qquad
    s=s_0+S(\mathcal L).
\]
These are $P_t^<(w)$ and $S_t^<(w)$, so $F_t(w)=w(p-\epsilon_0)-s$ is
evaluated in constant time. If it is nonnegative, retain the node and
its prefix totals as a candidate and descend left. Otherwise, replace
$(p_0,s_0)$ by $(p+p_w,s+s_w)$ and descend right. Empty subtrees have
zero aggregates. The last retained candidate is the smallest satisfying
key, and its saved totals give the threshold in constant additional time.
The descent visits $O(\log d)$ nodes.
\end{proof}

\begin{proof}[Proof of Lemma \ref{lem:upper-query}]
The crossing inequalities again ensure that $v$ exists and has a
predecessor. If $G_t(v)<0$, the root lies strictly between the predecessor
of $v$ and $v$. Throughout that interval the active upper set is exactly
$\{i:z_i(t)\geq v\}$, including the entire group at $v$. Hence
\[
    G_t(x)=S_t^{\geq}(v)-x\{P_t^{\geq}(v)+\epsilon_0\},
\]
and solving $G_t(x)=0$ gives the formula.

Use the same running-prefix descent as in the lower query. At a node
with key $w$ and prefix totals $p=P_t^<(w)$, $s=S_t^<(w)$, evaluate
\[
    G_t(w)=C_t-s-w(1-p+\epsilon_0).
\]
This identity follows because the suffix including the node at $w$ has
mass $1-p$ and count sum $C_t-s$, and the group at $w$ contributes zero
to $G_t(w)$. Retain the node and descend left when $G_t(w)\leq0$;
otherwise update the prefix totals and descend right. At the selected
key, the saved totals give
\[
    P_t^{\geq}(v)=1-P_t^<(v),
    \qquad
    S_t^{\geq}(v)=C_t-S_t^<(v).
\]
Thus each visited node requires constant work, and the entire query
takes $O(\log d)$ time.
\end{proof}

\begin{proof}[Proof of \Cref{thm:dynamic}]
We first specify the state maintained by the algorithm, verify that
initialization and each observation update preserve it, and then
deduce the claimed complexity and pathwise equality.

\medskip
\noindent
\emph{Maintained state and initialization.}
Immediately before observing $X_{t+1}$, the algorithm maintains the
following invariant:
\begin{enumerate}
    \item The category arrays and the stored total contain the exact
    values
    \[
        c_i(t)=N_i(t)+\alpha_i,
        \qquad
        z_i(t)=\frac{c_i(t)}{P_i},
        \qquad
        C_t=\sum_{i=1}^d c_i(t)=t+\alpha_+.
    \]

    \item The balanced tree has one node for each distinct key
    $z_i(t)$, with ties grouped together. All group totals,
    group cardinalities, and subtree aggregates are correct.

    \item The stored indicator correctly records whether
    $\widehat U_t\in\cP_0$. Whenever
    $\widehat U_t\notin\cP_0$, the stored values $A_t,B_t$
    are exactly its scaled oracle thresholds.
\end{enumerate}

At $t=0$, set $c_i(0)=\alpha_i$, $z_i(0)=\alpha_i/P_i$, and
$C_0=\alpha_+$. Sorting the initial keys, grouping ties, and building
the augmented balanced tree establishes the first two parts of the
invariant. Equation~\eqref{eq:dynamic-tv} determines null-ball
membership exactly. If the initial predictor is outside the null
ball, Lemmas~\ref{lem:lower-query} and~\ref{lem:upper-query} return
its exact scaled oracle thresholds. Thus the invariant holds at
$t=0$.

Forming the initial keys takes $O(d)$ time, sorting takes
$O(d\log d)$ time, and grouping ties and building the tree with its
subtree aggregates take $O(d)$ time. The membership query and,
when required, the two threshold queries take $O(\log d)$ time
altogether. Hence initialization costs $O(d\log d)$.

\medskip
\noindent
\emph{One observation update.}
Suppose the invariant holds at time $t$, and let $X_{t+1}=k$.
Using the current stored state, the algorithm first evaluates
\[
    e_{t+1}
    =
    \begin{cases}
        1,
        & \widehat U_t\in\cP_0,\\[3pt]
        \displaystyle
        \frac{\clip(z_k(t),A_t,B_t)}{C_t},
        & \widehat U_t\notin\cP_0.
    \end{cases}
\]
By the invariant and \eqref{eq:dynamic-factor},
$e_{t+1}=E_{\widehat U_t}(k)$. Reading $z_k(t)$ from the category
array and applying the stored thresholds require only a fixed
number of arithmetic operations and comparisons. Thus factor
evaluation costs $O(1)$. The algorithm updates the wealth by
\[
    M_{t+1}^{\mathrm{dyn}}
    =
    M_t^{\mathrm{dyn}}e_{t+1}
\]
\emph{before incorporating the new observation into the predictor}.

It then removes category $k$'s old contribution from the node at
$z_k(t)$, updates
\[
    c_k(t+1)=c_k(t)+1,
    \qquad
    z_k(t+1)=z_k(t)+\frac{1}{P_k},
    \qquad
    C_{t+1}=C_t+1,
\]
and inserts the updated contribution at $z_k(t+1)$. No other
category count or key changes. After repairing the affected group
totals and subtree aggregates, the arrays and tree therefore
represent exactly the predictor $\widehat U_{t+1}$.

Each removal or insertion takes $O(\log d)$ worst-case time.
Indeed, locating the relevant group and performing any necessary
node insertion, deletion, or rebalancing involve only
$O(\log d)$ nodes. Once a group is located, its totals and
cardinality are changed in constant time. At each affected node,
the subtree aggregates are recomputed from its group totals and
the two child-subtree aggregates, also in constant time. Thus
maintaining the augmentation does not increase the asymptotic
cost of the balanced-tree operations.

Finally, the algorithm applies \eqref{eq:dynamic-tv} to determine
membership for $\widehat U_{t+1}$. If this predictor is outside the
null ball, it uses Lemmas~\ref{lem:lower-query}
and~\ref{lem:upper-query} to obtain its exact scaled oracle
thresholds; otherwise it records that the next factor is
identically one. This establishes the invariant at time $t+1$.
These queries require at most three further tree descents, each
costing $O(\log d)$.

Consequently, a complete observation update consists of constant
work, one removal, one insertion, and a constant number of
logarithmic-time queries. Its total worst-case cost is therefore
$O(\log d)$.

\medskip
\noindent
The induction above shows that, on every round $t+1$, the dynamic
implementation applies exactly $E_{\widehat U_t}(X_{t+1})$, with
$\widehat U_t$ computed from the first $t$ observations. Therefore,
for every observation sequence,
\[
    M_n^{\mathrm{dyn}}
    =
    \prod_{t=0}^{n-1} E_{\widehat U_t}(X_{t+1})
    =
    M_n^{\mathrm{CK}},
\]
where the empty product equals one. This proves the claimed
identity for every $n\geq0$.

\medskip
\noindent
\emph{Storage and total cost.}
Each tree node represents a nonempty group of categories, so there
are at most $d$ nodes. Each node stores a constant number of
fields, and the category arrays contain $O(d)$ entries. The
global state requires constant storage, and any traversal stack
requires at most $O(\log d)$ additional words. Hence the total
storage is $O(d)$ words.

Combining the $O(d\log d)$ initialization cost with the
$O(\log d)$ worst-case cost of each of the $n$ observation updates
gives
\[
    O(d\log d+n\log d),
\]
as claimed.
\end{proof}

\begin{proof}[Proof of Proposition~\ref{prop:farrell-computation}]
\label{app:farrell-computation}
Maintain the global counts and the raw-count keys
\[
 v_i(n)=N_i(n)/P_i.
\]
Only the observed category changes its key, by $1/P_i$. As in the dynamic implementation in \Cref{sec:dynamic}, use an augmented balanced search tree, grouping equal keys. Store the group and subtree sums of $P_i$ and $N_i(n)$, together with group cardinalities. The all-zero initial keys form one group, so initialization costs $O(d)$. One removal and insertion maintain the tree in $O(\log d)$ worst-case time per observation.

At time $m=m_j$, one strict suffix query gives
\[
 p_j=\sum_{v_i(m)>m}P_i=P(\widehat A_j),\qquad
 s_j=\sum_{v_i(m)>m}N_i(m)=m\bar S_m(\widehat A_j).
\]
Consequently $q_j=p_j+\epsilon_0$ and
\[
 \widehat\Delta_j=(s_j/m-p_j-\epsilon_0)_+.
\]
No category list is needed for this query.

For membership during a betting block, maintain arrays $\mathrm{stamp}[i]$ and $\mathrm{bit}[i]$, initializing every stamp to $-1$. The first observation fills the first pilot, after which account zero is selected. The current block identifier is $j$. On observing category $i$ at time $t\in\{m_j+1,\ldots,2m_j\}$, \emph{before incrementing its count}, do the following:
\begin{enumerate}[itemsep=2pt,topsep=3pt]
 \item If $\mathrm{stamp}[i]\ne j$, set
 $\mathrm{bit}[i]=\one\{N_i(t-1)>m_jP_i\}$ and $\mathrm{stamp}[i]=j$.
 \item Use $\mathrm{bit}[i]$ to update the current account.
 \item Increment $N_i$ and repair its tree entry.
\end{enumerate}
At the first appearance of $i$ in the block, it has not been observed since the pilot ended. Therefore $N_i(t-1)=N_i(m_j)$, so the stored bit equals $\one\{i\in \widehat A_j\}$. Later appearances reuse that exact bit. At a new block, incrementing the identifier invalidates all old bits without clearing the arrays.

At a dyadic endpoint, first finish and freeze the old account, then select the next bet from the updated counts. One suffix query initializes the new account. Thus an observation requires constant account work, one logarithmic-time tree update, and at most one logarithmic-time query. Arrays and tree use $O(d)$ memory. Equation~\eqref{eq:farrell-finite-state} evaluates the mixture in constant time from the account state.

For the combined process, retain a separate tree for the smoothed keys used by clipped Kelly. Evaluate both betting factors from their pre-observation states, increment the shared counts once, and update both trees. At a dyadic endpoint, select the next learned event only after this update. The two trees have $O(d)$ total storage and each requires only a constant number of logarithmic-time operations per observation. Including the clipped-Kelly initialization in \Cref{thm:dynamic} gives total work $O(d\log d+n\log d)$ through time $n$. These counts include unit-cost exponential evaluations for the account factors. They make no bit-complexity or floating-point certification claim. 
\end{proof}

\subsection{Verification of the Farrell lower bound}
\label{app:farrell-lower}
Choose a category $i$ with $\epsilon_0<1-P_i$, and put $q=P_i+\epsilon_0\in(0,1)$. Define the full-support boundary null
\[
 R_i^0=q,\qquad R_\ell^0=P_\ell\frac{1-q}{1-P_i}\quad(\ell\ne i).
\]
For $0<\delta<1-q$, set
\[
 S_i^\delta=q+\delta,\qquad
 S_\ell^\delta=R_\ell^0\frac{1-q-\delta}{1-q}\quad(\ell\ne i).
\]
Then $\TV(R^0,P)=\epsilon_0$ and $\TV(S^\delta,P)=\epsilon_0+\delta$. Every null has $R_i\le q$, so data processing gives
\[
 \Gamma(S^\delta)\ge\inf_{r\le q}\kl(q+\delta,r)=\kl(q+\delta,q).
\]
The conditional distributions on $\cX\setminus\{i\}$ agree under $S^\delta$ and $R^0$. Therefore
\begin{equation}
 \Gamma(S^\delta)=\KL(S^\delta\|R^0)
 =\kl(q+\delta,q)\sim\frac{\delta^2}{2q(1-q)}.
 \label{eq:farrell-anchor}
\end{equation}
Choose a sufficiently small fixed $a>0$, and then $\delta_k\downarrow0$ with $\Gamma(S^{\delta_k})=ae^{-(2k+1)}$. This is a direct sequence with anchor $R^0$ and comparison constant one, in the terminology of~\cite[Section 3.2]{agrawal2025stopping}. That section explicitly permits fixed rescalings of the geometric buckets. Its lower bound implies a positive limsup of $\Gamma(S^{\delta_k})\E_{S^{\delta_k}}\tau/\ell(\Gamma(S^{\delta_k}))$ for every level-$\alpha$, power-one test, with $\alpha<1/2$. This proves \eqref{eq:farrell-lower}. The statement is sequence-wise; it supplies no universal lower bound at every isolated $S$. \hfill$\square$

\section{Mixture over pointwise oracles}
\label{app:mixture}
For completeness, we collect here the oracle-mixture construction and its regret guarantees.  These results are not needed for the main predictable plug-in procedure. This is an alternative method, which is computationally

Let $\mathcal A\subset\Delta_d\setminus\cP_0$ be compact, contained in
some $\cK_{m,\delta}$, and let $\Pi$ be a probability measure with full
support on $\mathcal A$.  Define the oracle mixture
\begin{equation}\label{eq:mixture}
    B_n^\Pi
    :=\int_{\mathcal A}\prod_{t=1}^n E_U(X_t)\,\Pi(du).
\end{equation}

\begin{theorem}[Finite-sample regret of the oracle mixture]\label{thm:mixture-regret}
The process $(B_n^\Pi)$ is a test supermartingale under the conditionally
varying null.  Fix a full-support $S\in\mathcal A$.  For every measurable
$V\subseteq\mathcal A$ with $\Pi(V)>0$,
\begin{equation}\label{eq:mixture-regret}
\begin{aligned}
    0
    &\leq n\Gamma(S)-\E_S\log B_n^\Pi\\
    &\leq \log\frac1{\Pi(V)}
      +n\sup_{U\in V}\{\Gamma(S)-h(S,U)\}\\
    &\leq \log\frac1{\Pi(V)}
      +n\sup_{U\in V}\KL(S\|U).
\end{aligned}
\end{equation}
\end{theorem}

\begin{proof}
Each fixed-$U$ product is a valid test supermartingale, so Tonelli's theorem
proves validity of the mixture. 
$B_n^\Pi$ is a valid horizon-$n$
e-variable under $R_S^n$. Set
\[
L_n=\prod_{t=1}^n\frac{S(X_t)}{R_S(X_t)}.
\]
Full support and Jensen’s inequality give
\[
\begin{aligned}
\mathbb E_S\log B_n^\Pi
&=n\Gamma(S)+
  \mathbb E_S\log\frac{B_n^\Pi}{L_n}\\
&\le n\Gamma(S)+
  \log\mathbb E_{R_S^{\otimes n}}B_n^\Pi\\
&\le n\Gamma(S).
\end{aligned}
\]

Let $\Pi_V$ be $\Pi$ conditioned on $V$.  Pathwise, Jensen's inequality for
the concave logarithm gives
\begin{align*}
    \log B_n^\Pi
    &\geq \log\Pi(V)
      +\log\int_V\prod_{t=1}^n E_U(X_t)\,\Pi_V(du)\\
    &\geq \log\Pi(V)
      +\int_V\sum_{t=1}^n\log E_U(X_t)\,\Pi_V(du).
\end{align*}
Taking expectation under $S$ yields
\[
    \E_S\log B_n^\Pi
    \geq
    \log\Pi(V)+n\int_V h(S,U)\,\Pi_V(du).
\]
Rearranging proves the second line of \eqref{eq:mixture-regret}; the last line
follows from Lemma~\ref{lem:oracle-mismatch}.
\end{proof}

\begin{corollary}[Logarithmic mixture regret under local prior mass]\label{cor:mixture-log-regret}
Fix a full-support $S\in\mathcal A$ and put $m_S=\min_i S_i$.  Suppose there
are constants $c_S>0$, $\kappa_S>0$, and $\rho_S>0$ such that
\begin{equation}\label{eq:small-ball}
    \Pi\!\left(\{U\in\mathcal A:\|U-S\|_1\leq\rho\}\right)
    \geq c_S\rho^{\kappa_S},
    \qquad 0<\rho\leq\rho_S.
\end{equation}
Then there is a finite constant $C_S$ such that, for all $n\geq2$,
\begin{equation}\label{eq:mixture-log-regret}
    0\leq n\Gamma(S)-\E_S\log B_n^\Pi
    \leq \frac{\kappa_S}{2}\log n+C_S.
\end{equation}
In particular,
\[
    \frac1n\E_S\log B_n^\Pi\longrightarrow\Gamma(S).
\]
For a prior with a density bounded below in a neighborhood of an interior
point $S$ of the $(d-1)$-dimensional simplex, one may take
$\kappa_S=d-1$.
\end{corollary}

\begin{proof}
For $\rho\leq m_S/2$ and $\|U-S\|_1\leq\rho$, every $U_i\geq m_S/2$.
Using $\KL(S\|U)\leq\sum_i(S_i-U_i)^2/U_i$ gives
\[
    \KL(S\|U)
    \leq \frac{2}{m_S}\|U-S\|_2^2
    \leq \frac{2}{m_S}\rho^2.
\]
Apply Theorem~\ref{thm:mixture-regret} to the $\ell_1$ ball of radius $\rho$
and use \eqref{eq:small-ball}.  Choosing
$\rho$ proportional to $n^{-1/2}$, truncated if necessary to lie below
$\min\{\rho_S,m_S/2\}$, yields \eqref{eq:mixture-log-regret}.
\end{proof}

\begin{corollary}[Almost-sure growth of the oracle mixture]\label{cor:mixture-as}
Under the assumptions above, for every fixed full-support $S\in\mathcal A$,
\[
    \frac1n\log B_n^\Pi\longrightarrow\Gamma(S)
    \qquad\text{almost surely under }S.
\]
\end{corollary}

\begin{proof}
Let $\widehat S_n$ be the empirical distribution.  Since each $E_U$ is
null-valid,
\[
    \frac1n\log B_n^\Pi
    \leq \Gamma(\widehat S_n)\longrightarrow\Gamma(S)
\]
almost surely.  Conversely, fix $\eta>0$.  Continuity provides a
neighborhood $V$ of $S$ on which $h(S,U)\geq\Gamma(S)-\eta$.  Uniform
continuity on the compact parameter set implies
$\inf_{U\in V}h(\widehat S_n,U)\geq\Gamma(S)-2\eta$ eventually.  Hence
$B_n^\Pi\geq\Pi(V)e^{n(\Gamma(S)-2\eta)}$ eventually.  Let
$\eta\downarrow0$.
\end{proof}




\subsection{Finite-grid approximation to the mixture}

Let $U^{(1)},\ldots,U^{(K)}$ be candidate alternatives and let
$\pi_k>0$, $\sum_k\pi_k=1$.  Then
\[
    B_n^{(K)}
    =\sum_{k=1}^K\pi_k\prod_{t=1}^n E_{U^{(k)}}(X_t)
\]
is a test supermartingale and, pathwise,
\[
    \log B_n^{(K)}
    \geq
    \max_{k\leq K}
    \left\{
      \sum_{t=1}^n\log E_{U^{(k)}}(X_t)+\log\pi_k
    \right\}.
\]
Thus a uniform prior pays at most $\log K$ relative to the best grid oracle.
Precomputing $K$ static oracles costs $O(Kd\log d)$ using sorting and each
observation costs $O(K)$ thereafter.  The statistical price is the
approximation gap
\[
    \Gamma(S)-\max_{k\leq K}\E_S\log E_{U^{(k)}}(X).
\]

A direct nonasymptotic stopping-time consequence is also available.  If
$\mu_k(S):=\E_S\log E_{U^{(k)}}(X)>0$ and
$L_k:=\max_i\log E_{U^{(k)}}(i)$, then the $k$th expert log-wealth is a random walk with bounded increments and positive drift $\mu_k(S)$, so its threshold-hitting time has finite expectation.  Wald's identity and the bounded overshoot then give
\[
    \E_S\tau_\alpha(B^{(K)})
    \leq
    \frac{
      \log(1/\alpha)+\log(1/\pi_k)+L_k
    }{\mu_k(S)}.
\]
Minimizing over $k$ gives the best finite-grid bound.

\end{document}